\documentclass[reqno,10pt,centertags,draft]{amsart}
\usepackage{amssymb,upref}

\newcommand{\bbN}{{\mathbb{N}}}
\newcommand{\bbR}{{\mathbb{R}}}

\newcommand{\bbZ}{{\mathbb{Z}}}
\newcommand{\bbC}{{\mathbb{C}}}

\newcommand{\calA}{{\mathcal A}}
\newcommand{\calB}{{\mathcal B}}
\newcommand{\calC}{{\mathcal C}}
\newcommand{\calD}{{\mathcal D}}

\newcommand{\calH}{{\mathcal H}}

\newcommand{\calK}{{\mathcal K}}

\newcommand{\calZ}{{\mathcal Z}}

\newcommand{\no}{\nonumber}
\newcommand{\lb}{\label}
\newcommand{\f}{\frac}

\newcommand{\ol}{\overline}

\newcommand{\what}{\widehat}

\newcommand{\spec}{\text{\rm{spec}}}
\newcommand{\rank}{\text{\rm{rank}}}
\newcommand{\ran}{\text{\rm{ran}}}
\newcommand{\ind}{\text{\rm{index}}}
\newcommand{\trind}{\text{\rm{trindex}}}
\newcommand{\gtr}{\text{\rm{gtr}}}
\newcommand{\sgn}{\text{\rm{sgn}}}
\newcommand{\dom}{\text{\rm{dom}}}
\newcommand{\esss}{\text{\rm{ess.spec}}}
\newcommand{\acs}{\text{\rm{ac.spec}}}

\newcommand{\bi}{\bibitem}

\newcommand{\beq}{\begin{equation}}
\newcommand{\eeq}{\end{equation}}
\newcommand{\ba}{\begin{align}}
\newcommand{\ea}{\end{align}}

\renewcommand{\Re}{\text{\rm Re}}
\renewcommand{\Im}{\text{\rm Im}}

\DeclareMathOperator{\tr}{tr}
\DeclareMathOperator*{\nlim}{n-lim}
\DeclareMathOperator*{\slim}{s-lim}

\allowdisplaybreaks
\numberwithin{equation}{section}

\newtheorem{theorem}{Theorem}[section]
\newtheorem{lemma}[theorem]{Lemma}
\newtheorem{corollary}[theorem]{Corollary}
\newtheorem{hypothesis}[theorem]{Hypothesis}
\theoremstyle{definition}
\newtheorem{definition}[theorem]{Definition}

\theoremstyle{remark}
\newtheorem{remark}[theorem]{Remark}

\begin{document}
\title[$\Xi$ Operator]{The $\Xi$ Operator and its Relation
to Krein's Spectral Shift Function}
\author[Gesztesy and Makarov]{Fritz Gesztesy and
Konstantin A.~Makarov}
\address{Department of Mathematics,
University of
Missouri, Columbia, MO
65211, USA}
\email{fritz@math.missouri.edu\newline
\indent{\it URL:}
http://www.math.missouri.edu/people/fgesztesy.html}
\address{Department of Mathematics, University of
Missouri, Columbia, MO
65211, USA}
\email{makarov@azure.math.missouri.edu\newline
\indent{\it URL:}
http://www.math.missouri.edu/people/kmakarov.html}

\begin{abstract}
We explore connections between Krein's spectral shift
function $\xi(\lambda,H_0,H)$ associated with the pair of
self-adjoint operators $(H_0,H)$, $H=H_0+V$ in a 
Hilbert space $\calH$ and the recently introduced concept of
a spectral shift operator $\Xi(J+K^*(H_0-\lambda-i0)^{-1}K)$
associated with the operator-valued Herglotz function
$J+K^*(H_0-z)^{-1}K$,
$\Im(z)>0$ in $\calH$, where $V=KJK^*$ and $J=\sgn(V)$. Our
principal results include a new representation for
$\xi(\lambda,H_0,H)$ in terms of an averaged index
for the Fredholm pair of self-adjoint spectral projections
$(E_{J+A(\lambda)+tB(\lambda)}((-\infty,0)),
E_J((-\infty,0)))$, $t\in\bbR$, where
$A(\lambda)=\Re(K^*(H_0-\lambda-i0)^{-1}K)$, $B(\lambda)=
\Im(K^*(H_0-\lambda-i0)^{-1}K)$ a.e.~Moreover, introducing
the new concept of a trindex for a pair of operators $(A,P)$
in $\calH$, where $A$ is bounded and $P$ is an
orthogonal projection, we prove that $\xi(\lambda,H_0,H)$
coincides with the trindex associated with the pair
$(\Xi(J+K^*(H_0-\lambda-i0)^{-1}K),\Xi(J))$.
 In addition, we discuss a
variant of the Birman-Krein formula relating
 the trindex of a pair of $\Xi$-operators and the Fredholm determinant
of the abstract scattering matrix.

We also provide  a generalization
of the classical Birman-Schwinger principle, replacing
the traditional eigenvalue counting functions by
appropriate spectral shift functions.
\end{abstract}

\maketitle

\section{Introduction} \lb{s1}
In order to facilitate a description of the content
of this paper we briefly introduce the notation
used throughout this manuscript. The open
complex upper half-plane is abbreviated by
$\bbC_+=\{z\in\bbC\,|\,\Im(z)>0\};$ the symbols $\calH,
\calK$ represent complex separable Hilbert spaces and
$I_\calH$ the corresponding identity operator in $\calH$.
Moreover, we denote by
$\calB(\calH)$,
$\calB_\infty(\calH)$, the Banach spaces of bounded and
compact operators in $\calH$ and by $\calB_p(\calH)$,
$p\geq 1$ the standard Schatten-von Neumann trace ideals
(cf., \cite{GK69}, \cite{Si79}). Real and
imaginary parts of an operator $T$ with
$\dom(T)=\dom(T^*)$ are
defined as usual by $\Re(T)=(T+T^*)/2$ and $\Im(T)=
(T-T^*)/(2i);$ the spectrum, essential and absolutely
continuous spectrum  of $T$ is abbreviated by
$\spec(T)$, $\esss(T)$, and $\acs(T)$, respectively.
For a self-adjoint operator $H=H^*$ in
$\calH$, the associated family of strongly
right-continuous orthogonal spectral
projections of $H$ in $\calH$ is denoted by
$\{E_H(\lambda)\}_{\lambda\in\bbR}$.

Before describing the content of each section, we briefly
summarize the principal results of this paper. Let $H_0$ and
$H=H_0+V$ be self-adjoint operators in $\calH$ with
$V=V^*\in\calB_1(\calH)$ and denote by $\xi(\lambda,H_0,H)$
Krein's spectral shift function associated with the pair
$(H_0,H)$, uniquely defined for a.e.~$\lambda\in\bbR$ by
$\xi(\cdot,H_0,H)\in L^1(\bbR;d\lambda)$ and
\begin{equation}
\tr((H-z)^{-1}-(H_0-z)^{-1})=-\int_\bbR d\lambda\,
\xi(\lambda,H_0,H)(\lambda -z)^{-2}, \quad
 z\in\bbC\backslash\bbR. \lb{1.1}
\end{equation}
Denoting by $\ind (P,Q)$ the index of a Fredholm pair of
orthogonal projections in $\calH$, one of our principal
results represents $\xi(\lambda,H_0,H)$ as an averaged
 index for the Fredholm pair of projections
$(E_{J+A(\lambda)+tB(\lambda)}((-\infty,0)),
E_J((-\infty,0)))$, $t\in\bbR$, as follows (cf.
Theorem~\ref{main}),
\begin{align}
&\xi(\lambda, H_0,H)=\frac{1}{\pi} \int_\bbR dt
\frac{\ind \big(E_{J+A(\lambda)+tB(\lambda)}
((-\infty,0)),E_J((-\infty,0))\big )}{1+t^2}
\lb{1.2} 
\end{align}
for  a.e.~$\lambda\in\bbR$. Here $V=H-H_0$ is decomposed as
\begin{equation}
V=KJK^*, \quad K\in\calB_2(\calH), \,\, J=\sgn(V) \lb{1.3}
\end{equation}
and
\begin{align}
&A(\lambda)=\nlim_{\varepsilon \downarrow 0}
 \Re(K^*(H_0-\lambda-i\varepsilon)^{-1}K), \lb{1.4} \\
&B(\lambda)=\nlim_{\varepsilon \downarrow 0}
\Im(K^*(H_0-\lambda-i\varepsilon)^{-1}K)
 \text{ for a.e. } \lambda\in\bbR. \lb{1.5}
\end{align}
In the special sign-definite case where $J=\pm I_\calH$, formula
\eqref{1.2} implies a result of Pushnitskii \cite{Pu97}
(cf. Corollary~\ref{c5.6}).

Next, observing that the trace class-valued operator
$K^*(H_0-z)^{-1}K$, $z\in\bbC_+$ has nontangential boundary
values $K^*(H_0-\lambda-i0)^{-1}K$ for a.e.~$\lambda\in\bbR$
in $\calB_p(\calH)$-topology for each $p>1$, but in general
not in the trace norm $\calB_1(\calH)$-topology, we
introduce the notion of a trindex, $\trind (\cdot,\cdot)$,
for a pair of operators $(A,Q)$ in $\calH$, where
$A\in\calB(\calH)$ is bounded and $Q=Q^*=Q^2$ is an
orthogonal projection in $\calH$, as follows: the pair
$(A,Q)$ is said to have a {\it trindex}, denoted by $\trind
(A,Q)$, if there exists an orthogonal projection $P$ in
$\calH$ such that $(A-P)\in\calB_1(\calH)$ and $(P,Q)$ is a
Fredholm pair of orthogonal projections in $\calH$. In this
case one then defines,
\begin{equation}
\trind (A,Q)=\tr(A-P)+\ind (P,Q). \lb{1.6}
\end{equation}
Introducing the spectral shift operator
$\Xi(T)=(1/\pi)\Im(\log(T))$ for a bounded dissipative
operator $T$, with $T^{-1}$ also bounded in $\calH$, our second
principal result (cf. Theorem~\ref{t5.8}) identifies
$\xi(\lambda,H_0,H)$ with the trindex of the pair
$(\Xi(J+K^*(H_0-\lambda-i0)^{-1}K),\Xi(J))$, that is,
\begin{equation}
\xi(\lambda,H_0,H)=\trind(\Xi(J+K^*(H_0-
\lambda-i0)^{-1}K), \Xi(J)) \text{ for a.e. }
\lambda \in\bbR. \lb{1.7}
\end{equation}
The trindex representation \eqref{1.7} paves the way for 
introducing a generalized
spectral shift function to be discussed in
detail in Section~\ref{s5}. We also show that
the averaging formula
\eqref{1.2} should be viewed as a generalized Birman-Schwinger
principle
\begin{equation}
\xi(\lambda, H_0,H)= \frac{1}{\pi}\int_\bbR dt\,
\frac{\what\xi(0_-,J+A(\lambda)+tB(\lambda),J)}{1+t^2}
 \text{ for a.e. } \lambda \in \bbR, \lb{1.5a}
\end{equation}
where $\what\xi(\cdot,J+A(\lambda)+tB(\lambda),J)$
denotes the
generalized spectral shift function associated with the pair
$(J+A(\lambda)+tB(\lambda),J)$,  $t\in \bbR$.

Next we turn to a description of the content of each
section. In Section~\ref{s2} we briefly review basic
properties of the index of a Fredholm pair of projections in
$\calH$ (mainly following \cite{ASS94}, see also
\cite{AS94}, \cite{Ka97}, \cite{Ka55}) and then present a
discussion of the notion of trindex and some of its
properties. Section~\ref{s3} is devoted to the concept of a
$\Xi$-operator as recently introduced in \cite{GMN99} and
further discussed in \cite{GM99}. More precisely, if $T$ is
a bounded dissipative operator with $T^{-1}\in\calB(\calH)$,
then $\Xi(T)$ is defined by
\begin{equation}\lb{1.8}
\Xi(T)=\frac{1}{\pi} \Im (\log (T)).
\end{equation}
Section~\ref{s3} also studies sufficient conditions on
$A=A^*\in\calB_\infty(\calH)$ and $0\leq B\in
\calB_1(\calH)$ to
guarantee $(\Xi(S+A+iB)-\Xi(S+A))\in\calB_1(\calH)$ for
given $S=S^*\in\calB(\calH)$. Section~\ref{s4}, the
technical core of this paper, provides a discussion of
averaged Fredholm indices. Introducing the family of normal
trace class operators,
\begin{equation}
\calA (z)=B^{1/2}(S+zB)^{-1}B^{1/2}, \quad
z\in\bbC\backslash\bbR, \lb{1.9}
\end{equation}
where $S=S^*\in\calB(\calH)$, $S^{-1}\in\calB(\calH)$, $0\leq
B\in\calB_1(\calH)$,
associated with the (dissipative) family of operator-valued
Herglotz functions
\begin{equation}
T(z)=S+zB, \quad z\in\bbC_+ \lb{1.10}
\end{equation}
(i.e., $T$ is analytic in $\bbC_+$ and $\Im(T(z))\geq 0$
for $z\in\bbC_+$), we prove
\begin{equation}\lb{1.11}
\tr (\log(T(z))-\log(S))=
\sum_{k=1}^\infty m_k \int_0^1 d\tau\, \frac{z
\lambda_k}{1+\tau z \lambda_k},
\end{equation}
(cf.~Theorem~\ref{tt.2}), where $\{\lambda_k \}_{k\in\bbN}
\subset \bbR$ is the set of eigenvalues with associated
multiplicities $\{m_k\}_{k\in\bbN}$ of the self-adjoint trace
class operator $B^{1/2}S^{-1}B^{1/2}$. This then yields our
principal result (Theorem~\ref{ttr.8}) relating the trindex
of the pair  $(\Xi(S+A+iB), \Xi(S))$ and the averaged Fredholm
index $n(t)=\ind (\Xi(S+A+tB),\Xi(S))$, $t\in\bbR$, as
\begin{equation}\lb{1.12}
\trind
(\Xi( S+A+iB),
\Xi(S))=
\frac{1}{\pi}\int_\bbR \frac{dt\,n(t)}{1+t^2}.
\end{equation}
Moreover, we provide a version of the Birman-Krein
 formula relating the left-hand side of
\eqref{1.12} and the determinant of the abstract scattering matrix
\begin{equation}
\exp (-2\pi i\, \trind
(\Xi( S+A+iB),
\Xi(S)))=\det (I_\calH-2iB^{1/2}(S+A+iB)^{-1}B^{1/2}).
\end{equation}
Section~\ref{s5} finally presents the applications of our
formalism to Krein's spectral shift function
$\xi(\lambda,H_0,H)$ as discussed in \eqref{1.2},
\eqref{1.7}, and \eqref{1.5a}.


In conclusion, we note that Krein's spectral shift function
\cite{Kr53}--\cite{KJ81}, a concept originally
introduced by
Lifshits \cite{Li52}, \cite{Li56}, continues to generate
considerable interest. Without repeating the extensive
bibliography recently provided in \cite{GMN99}, we remark
that the spectral shift function plays a fundamental role in
scattering theory, (relative) index theory, spectral
averaging and its application to localization properties of
random Hamiltonians, eigenvalue counting functions and
spectral asymptotics, semi-classical approximations, and
trace formulas for one-dimensional Schr\"odinger and Jacobi
operators. A very selective list of recent pertinent
references includes, for instance, \cite{BP98},
\cite{CHM96}, \cite{EP97}--\cite{GM99},
\cite{GS96}, \cite{KS98}, \cite{Mu98},
 \cite{Pu97}--\cite{Pu98a}, \cite{Ro98}, \cite{Sa97},
\cite{Si95},
\cite{Si98}. For many more references the interested
reader can
consult the 1993 reviews by Birman and Yafaev \cite{BY93},
\cite{BY93a}, and
 \cite{GM99},
\cite{GMN99}.

\section{Index of a Pair of Projections and Trindex}
\lb{s2}

In this section we recall the main properties of the
index of a Fredholm pair of orthogonal projections in
$\calH$ and discuss a closely related notion of a trindex
of a
pair  of operators, one of which is a bounded
operator in
$\calH$ and the other is an orthogonal projection in
$\calH$. Finally, we introduce the notion of a generalized
trace for a pair of bounded operators in $\calH$.

Let $P$ and $Q$ be orthogonal projections in a
complex separable  Hilbert space $\calH$. The pair
$(P,Q)$ is
said to be a  Fredholm pair if the map
\begin{equation}
QP\big|_{\ran(P)}:\ran(P) \to \ran (Q) \lb{2.1}
\end{equation}
is a Fredholm operator from the Hilbert space
 $\ran (P)$ to the Hilbert space $\ran (Q)$. In this case
one defines the index of the pair
$(P,Q)$ as the Fredholm
index of the operator $QP\big|_{\ran(P)}$
\begin{equation}\lb{tr.1}
\ind (P,Q)=\ind(QP\big|_{\ran(P)}).
\end{equation}

The following three results, Lemmas~\ref{l2.1} and \ref{ltr.4} 
and Theorem~\ref{ttr.3}, recall well-known
results for Fredholm pairs of projections. We refer, for
instance, to \cite{AS94}, \cite{ASS94}, \cite{Ka97},
\cite{Ka55} and the references cited therein.

We start with two important criteria for a pair $(P,Q)$ of
self-adjoint projections to be a Fredholm pair.

\begin{lemma}\lb{l2.1}
(i) A necessary and sufficient condition that $(P,Q)$ be
a Fredholm pair is that $P-Q=F+D$,
where  $F,D$ are self-adjoint,
$\|D\|<1$, and $F$ is a finite-rank operator.\\
(ii) $(P,Q)$ is a Fredholm pair if and only if  $+1$
and $-1$
do not belong to the essential spectrum of $(P-Q)$
\begin{equation}\lb{tr.2}
\pm 1\notin\esss (P-Q).
\end{equation}
In this case $\ker(P-Q\pm I_\calH)$ are both finite
dimensional and
\begin{equation}\lb{tr.3}
\ind (P,Q)=\dim (\ker (P-Q-I_\calH))-
\dim (\ker (P-Q+I_\calH)).
\end{equation}
In particular, if either
\begin{equation}\lb{tr.4}
(P-Q)\in \calB_\infty(\calH),
\end{equation}
or 
\begin{equation}\lb{tr.5}
\|P-Q\|<1,
\end{equation}
then $(P,Q)$ is a Fredholm pair.
\end{lemma}

The following result summarizes some of the most important
properties of the index of a Fredholm pair of projections.
\begin{theorem}\lb{ttr.3}
 (i) Let $(P,Q)$ be a Fredholm pair of projections in
$\calH$.
Then so is $(Q,P)$ and
\begin{equation}\lb{tr.6}
\ind (P,Q)=-\ind (Q,P).
\end{equation}
(ii) Let $(P,Q)$ and $(Q,R)$ be Fredholm pairs in
$\calH$ and
either $(P-Q)\in \calB_\infty(\calH)$ or $(Q-R)\in
\calB_\infty(\calH)$.  Then $(P,R)$ is a Fredholm pair
and one
has the chain rule
\begin{equation}\lb{tr.7}
\ind (P,R)=\ind (P,Q)+\ind (Q,R).
\end{equation}
(iii) If $\|P-Q\|<1$ then
\begin{equation}\lb{tr.8}
\ind (P,Q)=0.
\end{equation}
(iv) If $(P-Q)\in \calB_1(\calH)$ then
\begin{equation}\lb{tr.9}
\ind (P,Q)=\tr (P-Q).
\end{equation}
\end{theorem}

As shown in \cite{ASS94}, the compactness assumption in
connection with \eqref{tr.7} cannot be dropped in general.

Theorem~\ref{ttr.3}\,(ii)
combined with Lemma~\ref{l2.1}\,(i) implies a
stability result for the index for a Fredholm pair of
 projections under small perturbations.
\begin{lemma}\lb{ltr.4}
Let $P, $ $P_1$, and $Q$ be orthogonal projections
in $\calH$.
Assume that
\begin{equation}\lb{tr.10}
(P-Q)\in \calB_\infty(\calH)
\end{equation}
and
\begin{equation}\lb{tr.11}
\|P-P_1\|<1.
\end{equation}
Then $(P,Q)$ and $(P_1, Q)$ are Fredholm pairs in
$\calH$ and
\begin{equation}\lb{tr.12}
\ind (P,Q)=\ind (P_1,Q).
\end{equation}
\end{lemma}
\begin{proof}
That $(P,Q)$ is a Fredholm pair follows from
\eqref{tr.10} and Lemma~\ref{l2.1}\,(ii). Since
\eqref{tr.10} and \eqref{tr.11} hold, the difference
$P_1-Q$
can be represented as a sum of a contraction (with norm 
strictly less than one)
and a finite-rank operator. Thus 
 $(P_1,Q)$ is a Fredholm pair by Lemma~\ref{l2.1}\,(i). 
Since \eqref{tr.10}
holds one
can apply  Theorem~\ref{ttr.3}\,(ii) to
conclude
that
\begin{equation}\lb{tr.13}
\ind (P,P_1)=\ind (P,Q)+\ind (Q,P_1).
\end{equation}
By Theorem~\ref{ttr.3}\,(iii) and \eqref{tr.11}
one gets
\begin{equation}\lb{tr.14}
\ind (P,P_1)=0.
\end{equation}
Combining \eqref{tr.6}, \eqref{tr.13}, and
\eqref{tr.14} one arrives at \eqref{tr.12}.
\end{proof}

\begin{definition} \lb{d2.4} Let $A\in \calB(\calH)$
and $Q$ be
an orthogonal projection in $\calH$.
We say that the pair $(A,Q)$ has a {\it trindex}, denoted
by $\trind(A,Q)$, if there
exists an orthogonal projection $P$ in $\calH$ such that
 $(A-P)\in \calB_1(\calH)$ and $(P,Q)$ is a Fredholm pair
of orthogonal projections in
$\calH$. In this case the trindex of the pair $(A,Q)$ is
defined by
\begin{equation}\lb{tr.17}
\trind (A,Q)=\tr(A-P)+\ind (P,Q).
\end{equation}
\end{definition}

The following result shows that the trindex of the pair
$(A,Q)$ is well-defined, that is, it is independent of
the choice of the projection $P$
satisfying the conditions in Definition~\ref{d2.4}.

\begin{lemma}\lb{ltr.5}
Let  $A\in \calB(\calH)$, and  $P_1$, $P_2$,  and  $Q$
be orthogonal projections in $\calH$ such that
\begin{equation}\lb{tr.18}
(A-P_j) \in \calB_1(\calH), \quad j=1,2,
\end{equation}
and
\begin{equation}\lb{tr.19}
(P_j,Q) \text{ is a Fredholm pair, } \quad j=1,2.
\end{equation}
Then
\begin{equation}\lb{tr.20}
\tr(A-P_1)+\ind(P_1,Q)=\tr(A-P_2)+\ind(P_2,Q).
\end{equation}
\end{lemma}
\begin{proof}
By \eqref{tr.18} one concludes that $(P_1-P_2)\in
\calB_1(\calH)$ and hence by Lemma~\ref{l2.1}\,(ii),
the pair $(P_1,P_2)$ is a Fredholm pair. In particular,
Theorem~\ref{ttr.3}\,(iv) implies
\begin{equation}
\tr(A-P_1)=\tr(A-P_2)+\tr(P_2-P_1)
=\tr(A-P_2)+\ind(P_2,P_1).\lb{tr.21}
\end{equation}
Hence, Theorem~\ref{ttr.3}\,(ii) yields
\begin{align}
&\tr(A-P_1)+\ind(P_1,Q)
=\tr(A-P_2)+\ind(P_2,P_1)+\ind(P_1,Q) \no \\
&=\tr(A-P_2)+\ind(P_2,Q)\lb{tr.22}
\end{align}
proving \eqref{tr.20}.
\end{proof}

\begin{remark} \lb{r2.5}
Our motivation for introducing the concept of a trindex
for
a pair $(A,Q)$ lies in the following two facts.\\
(i) If $A=P$, with $(P,Q)$ a Fredholm pair of
projections in
$\calH$, then
\begin{equation}
\trind(P,Q)=\ind(P,Q). \lb{2.19}
\end{equation}
(ii) If $A\in\calB(\calH)$ and $Q$ is an
orthogonal projection in $\calH$ with
$(A-Q)\in\calB_1(\calH)$, then
\begin{equation}
\trind(A,Q)=\tr(A-Q). \lb{2.20}
\end{equation}
\end{remark}

The stability result for the trindex of a pair $(A,Q)$
analogous to Lemma~\ref{ltr.4} then reads as follows.

\begin{lemma} \lb{l2.6}
Let  $A\in \calB(\calH)$, and  $Q$, $Q_1$ be orthogonal
projections in $\calH$ such that $(Q-Q_1)\in
\calB_\infty(\calH)$ and
\begin{equation}
\|Q-Q_1\|<1.
\end{equation}
If $(A,Q)$ has a trindex, then $(A,Q_1)$ has a trindex and
\begin{equation}
\trind(A,Q)=\trind(A,Q_1).
\end{equation}
\end{lemma}
\begin{proof}
It suffices to combine Lemma~\ref{ltr.4} and \eqref{tr.17}.
\end{proof}

\begin{definition} \lb{sp.4} Let $A, B\in \calB(\calH)$
and $(A-B)\in\calB_\infty(\calH)$.
We say that the pair $(A,B)$ has a {\it generalized trace,}
denoted by $\gtr(A,B)$, if there
exists an orthogonal projection $Q$ in $\calH$ such that
both pairs, $(A,Q)$ and $(B,Q)$ have a trindex.
 In this case the generalized trace of the pair
$(A,B)$ is defined by
\begin{equation}\lb{sp17}
\gtr(A,B)=\trind(A,Q)-\trind (B,Q).
\end{equation}
\end{definition}

The following result shows that $\gtr(A,B)$ is well-defined,
that is,
it is independent of the choice of the orthogonal
projection $Q$
satisfying the conditions in Definition~\ref{sp.4}.

\begin{lemma} \lb{l2.9}
Let $A,B\in\calB(\calH)$, $(A-B)\in\calB_\infty(\calH)$,
and $Q_1,Q_2$ orthogonal projections in $\calH$ such that
$(A,Q_j)$ and $(B,Q_j)$, $j=1,2$ have a trindex. Then
\begin{equation}
\trind(A,Q_1)-\trind(B,Q_1)=\trind(A,Q_2)-\trind(B,Q_2).
\lb{2.27}
\end{equation}
\end{lemma}
\begin{proof}
By hypothesis, there exist orthogonal projections $P_{j,A},
P_{j,B}$ in $\calH$ such that $(A-P_{j,A}),(B-P_{j,B})\in
\calB_1(\calH)$ and the pairs $(P_{j,A},Q_j)$ and
$(P_{j,B},Q_j)$, $j=1,2$ are Fredholm pairs of projections.
Since $(A-B)\in\calB_\infty(\calH)$ one infers $(P_{j,A}
-P_{k,B})\in\calB_\infty(\calH)$, $j,k\in\{1,2\}$ and hence
$(P_{j,A},P_{k,B})$, $j,k\in\{1,2\}$ are Fredholm pairs.
Moreover, $(P_{2,A}-P_{1,A}),(P_{2,B}-P_{1,B})
\in\calB_1(\calH)$. Thus,
\begin{align}
&\trind(A,Q_1)-\trind(B,Q_1)-(\trind(A,Q_2)-\trind(B,Q_2))
\no \\
&=\tr(A-P_{1,A})+\ind(P_{1,A},Q_1)-\tr(A-P_{2,A})
-\ind(P_{2,A},Q_2) \no \\
&-\tr(B-P_{1,B})-\ind(P_{1,B},Q_1)+\tr(B-P_{2,B})
+\ind(P_{2,B},Q_2) \no \\
&=\tr((A-B-P_{1,A}+P_{1,B})-(A-B-P_{2,A}+P_{2,B}))
+\ind(P_{1,A},Q_1) \no \\
&-\ind(P_{1,B},Q_1) +\ind(P_{2,B},Q_2)
-\ind(P_{2,A},Q_2) \no \\
&=\tr(P_{2,A}-P_{1,A})-\tr(P_{2,B}-P_{1,B}) +
\ind(P_{1,A},P_{1,B})+\ind(P_{2,B},P_{2,A}) \no \\
&=\ind(P_{2,A},P_{1,A}) +\ind(P_{1,A},P_{1,B})
-\ind(P_{2,B},P_{1,B})  +\ind(P_{2,B},P_{2,A}) \no \\
&=\ind(P_{2,A},P_{1,B}) +\ind(P_{1,B},P_{2,A}) =0 \lb{2.28}
\end{align}
by repeatedly using \eqref{tr.6} and \eqref{tr.7}.
\end{proof}

\begin{remark} \lb{r2.10}
If $A,B\in\calB(\calH)$ with $(A-B)\in\calB_1(\calH)$,
then
\begin{equation}
\gtr(A,B)=\tr(A-B). \lb{2.29}
\end{equation}
\end{remark}

We were somewhat hesitant to introduce concepts
such as {\it trindex}, $\trind(\cdot,\cdot)$ and {\it
generalized trace},
$\gtr(\cdot,\cdot)$ as additional entities to such familiar
quantities like the Fredholm index, $\ind(\cdot,\cdot)$ and
the trace, $\tr(\cdot)$. However, it will become clear from
the remainder of this paper, that both objects seem to be
very natural in the context of  Krein's spectral shift
function (cf. Corollary~\ref{iff} and Remark~\ref{rindex}).

\section{The $\Xi$ Operator} \lb{s3}

 Suppose $T$ is a bounded dissipative operator in the
Hilbert
space $\calH$ (i.e., $\Im(T)\geq 0)$ and
 $L$ is the minimal self-adjoint dilation
of $T$ (cf.~\cite[Ch.~III]{SF70}) in the Hilbert space
$\calK\supseteq
\calH$.
 We define the $\Xi$-operator associated with the dissipative
operator $T$ by
\begin{equation} \lb{x.1}
\Xi(T)= P_{\calH} E_L((
-\infty,0))P_{\calH}\vert_{\calH},
\end{equation}
where $P_{\calH}$ is the orthogonal projection in
$\calK$ onto $\calH$ and
$\{ E_L(\lambda)\}_{\lambda\in \bbR}$ represents the
family of strongly right-continuous
orthogonal spectral projections of $L$ in $\calK$.

In particular, if $T=T^*$, the $\Xi$-operator coincides
with the spectral projection of $T$
corresponding to the  negative semi-axis $(-\infty, 0)$,
\begin{equation}\lb{x.3}
\Xi(T)= E_T(( -\infty,0)),
\end{equation}
since in this case  the minimal self-adjoint dilation of
$T$ coincides with $T$.

\begin{remark} \lb{r3.1}
(i) By \eqref{x.1}, $\Xi$ is a nonnegative contraction,
\begin{equation}\lb{x.4}
0\le \Xi(T)\le I_{\calH}.
\end{equation}
(ii) If $T\in \calB(\calH)$ is a bounded
dissipative operator and $T^{-1}\in \calB(\calH)$
then $\Xi(T)$ can be expressed in terms of the operator
logarithm of $T$ by
\begin{equation}\lb{x.2}
\Xi(T)=\pi^{-1} \Im (\log (T)),
\end{equation}
where $\log(T)$ is defined by
\begin{equation}\lb{2.17}
\log (T)=-i\int_0^\infty  d \lambda \,
((T+i\lambda)^{-1}-(1 +i\lambda)^{-1}I_\calH)
\end{equation}
in the sense of a $\calB(\calH)$-norm convergent Riemann
integral (cf. the extensive treatment in \cite{GMN99}).
Without going into details (these may be found in
\cite{GMN99}), we remark that \eqref{x.2} resembles the
exponential Herglotz representation
for scalar-valued Herglotz functions
studied in detail by
Aronszajn and Donoghue \cite{AD56}.
 Indeed, for any Herglotz
function $t(z)$ (i.e., $t:\bbC_+\to\bbC_+$ analytically),
$\log(t(z))$ is also a Herglotz function admitting the
representation,
\begin{equation}
\log(t(z))=c+\int_{{\mathbb{R}}}d\lambda \,
\xi(\lambda)((\lambda-z)^{-1}-\lambda
(1+\lambda^2)^{-1}), \quad z\in\bbC_+, \lb{2.18}
\end{equation}
where $c\in \bbR$ and
\begin{equation}\lb{2.19a}
0\leq \xi \leq 1 \text{ and } \xi
(\lambda)=\pi^{-1}\lim_{\varepsilon
\downarrow 0}\Im (\log (t(\lambda+i\varepsilon)))
\text{ a.e.}
\end{equation}
\end{remark}

At this point a natural question arises. Suppose $S$ is
a bounded dissipative
operator in $\calH$ and $T=S+A+iB$, $A=A^*\in\calB(\calH)$,
$0\leq B\in\calB(\calH)$ its dissipative bounded
perturbation. Can one expect an interesting relationship
between $\Xi(T)$ and $\Xi(S)$? The following is a first
result in this direction.

\begin{lemma}\lb{lx.0}
Let  $S,T \in \calB(\calH)$ be dissipative operators such
that
$S^{-1}, T^{-1} \in \calB(\calH)$ and assume
$(T-S)\in \calB_1(\calH)$. Then
\begin{equation}\lb{x.4b}
(\Xi(T) -\Xi(S))\in \calB_1(\calH).
\end{equation}
\end{lemma}
\begin{proof}
Since $S$ and $T$ have a bounded inverse, $\log(T)$ and
$\log(S)$ are well-defined and (cf.~\cite{GMN99})
\begin{align}
&\log(T)-\log(S)=
-i\int_0^\infty dt\,
 ((T+it)^{-1}-(S+it)^{-1}) \no \\
&=i\int_0^\infty dt\, ((T+it)^{-1}(T-S)(S+it)^{-1}).
\lb{x.4c}
\end{align}
Using standard estimates for resolvents
$(T+it)^{-1}$ and $(S+it)^{-1}$ $(t\ge 0)$
of the dissipative operators $T$ and $S$ entering
\eqref{x.4c} (cf.~Lemma~2.6 in \cite{GMN99})
one concludes
\begin{equation}\lb{x.4d}
(\log(T)-\log(S) )\in \calB_1(\calH),
\end{equation}
if $(T-S)\in \calB_1(\calH)$. Thus,
$(\Xi(T)-\Xi(S)) \in \calB_1(\calH)$ by \eqref{x.2} and
\eqref{x.4d}.
\end{proof}

We assume the following hypothesis in the sequel.

\begin{hypothesis}\lb{hx.1}
Assume $S=S^*\in\calB(\calH)$,
$A=A^*\in\calB_\infty(\calH)$,
 $0\leq B\in \calB_1(\calH)$,  and
$(S+A+\tau_0 B)^{-1}\in\calB(\calH)$ for some
$\tau_0\in \bbR$.
\end{hypothesis}

\begin{remark}\lb{rx.1}
Under Hypothesis~\ref{hx.1} one infers
\begin{equation}\lb{ess}
0\notin \esss(S)
\end{equation}
by the stability of the essential spectrum under compact
perturbations. Moreover, by the analytic Fredholm
theorem (cf.~\cite[Sect.~VI.5]{RS80}),
$(S+A+zB)^{-1}\in\calB(\calH)$,
$z\in\overline {\bbC_+}$  except for $z$ in a discrete set
$\calD\subset\bbR$, with $\pm\infty$ the only
possible accumulation points of $\calD$. In particular,
\begin{equation}\lb{bound}
(S+A+\varepsilon B)^{-1} \in \calB(\calH) \text{ for }
\varepsilon \text{ sufficiently small, } \varepsilon \ne 0.
\end{equation}

\end{remark}

\begin{lemma}\lb{lx.1}
Assume Hypothesis~\ref{hx.1} with $A=0$. Then
\begin{equation}\lb{x.5}
(\Xi(S+t B)-\Xi(S))\in \calB_1(\calH) \text{ for all }
t \in
\bbR
\end{equation}
and
\begin{equation}\lb{x.6}
\|\Xi(S+t B) -
\Xi(S)\|_{\calB_1(\calH)}=O(t) \text{ as } t \downarrow 0.
\end{equation}
Moreover, $S+zB$, $z\in\bbC_+$ has a bounded inverse,
\begin{equation}\lb{x.6a}
(S+zB)^{-1} \in \calB(\calH) \text{ for all } z\in\bbC_+
\end{equation}
and
\begin{equation}\lb{x.6b}
(\Xi(S+z B)-\Xi(S))\in \calB_1(\calH) \text{ for all } z\in
\ol {\bbC_+}\,.
\end{equation}
\end{lemma}
\begin{proof} Since $(S+\tau_0 B)^{-1}\in\calB(\calH)$ for
some $\tau_0\in \bbR$ and $B\in\calB_\infty (\calH)$, one
infers that
$0\notin \esss (S)$. Thus,
there exists a closed interval
$\Delta$, $0\in \Delta$, such that
 $\dim(\ran (E_{S+t B}(\Delta)))<\infty$ for all $t\in \bbR$.
 Hence,
for given
$t\in\bbR$, one can find a clockwise oriented bounded contour
$\Gamma_t$
encircling  $(\spec(S+t B)\cup\spec(S))\cap
(-\infty,0)$ such that
\begin{equation}\lb{x.7}
\Xi(S+t B)=E_{S+t B}((-\infty,0))=\frac{1}{2\pi i}
\oint_{\Gamma_t}
d\zeta\, (S+t B-\zeta)^{-1}, \quad t\in\bbR
\end{equation}
and
\begin{equation}\lb{x.8}
\Xi(S)=E_{S}((-\infty,0))=\frac{1}{2\pi i}\oint_{\Gamma_t}
d\zeta\, (S-\zeta)^{-1}.
\end{equation}
The second resolvent identity for $S+t B$ and $S$ then
implies
 $((S+tB-\zeta)^{-1}-(S-\zeta)^{-1})\in \calB_1(\calH)$,
$\zeta\in \Gamma_t$, proving \eqref{x.5}.

Since $B\geq 0$, there exists a closed interval
$\Delta$, $ \Delta\subset (-\infty,0)$, such that
\begin{equation}\lb{x.9}
\bigcup_{t\in [0, \varepsilon]}\spec
(S+t B) \cap (-\infty,0) \subset \Delta \text{ for }
\varepsilon >0 \text{ sufficiently small.}
\end{equation}
Thus, choosing the contour
\begin{equation}\lb{x.10}
\Gamma=\{\zeta\in \bbC \, | \, \text{dist}
(\zeta, \Delta)=\frac{1}{2} \text{dist}(0,\Delta)\},
\end{equation}
the representation \eqref{x.7} is valid
for all $t\in [0,\varepsilon]$, for sufficiently small
$\varepsilon >0$. Using the second resolvent
identity again and the standard estimate
$\|(A-\zeta)^{-1}\|\le (\text{dist}(\zeta,
\spec(A)))^{-1}$
 for every self-adjoint operator $A$ in $\calH$, one
obtains
\begin{align}
\|\Xi(S+t B)-\Xi(S) \|_{\calB_1(\calH)}&=
\|E_{S+t B}((-\infty, 0))
-E_{S}((-\infty, 0)) \|_{\calB_1(\calH)}
\no \\
&
\le (2/\pi)t (\text{dist}(0,\Delta))^{-2}
\| B\|_{\calB_1(\calH)} |\Gamma|,
\quad t \in [0,\varepsilon], \lb{x.11}
\end{align}
where $|\Gamma|$ denotes
the length of the contour $\Gamma$. This proves
\eqref{x.6}.

Next, consider the operator-valued Herglotz function
\begin{equation}\lb{x.12}
M(z)=S+zB, \quad z\in\bbC_+
\end{equation}
(i.e., $\Im(M(z))\geq 0$ for all $z\in\bbC_+$). By
hypothesis
there exists a $\tau_0\in\bbR$ such that
$(S+\tau_0 B)^{-1}\in\calB(\calH)$. Thus, for sufficiently
small values of $|z-\tau_0|$, $z\in \bbC_+$, the operator
$M(z)$,
being a small perturbation of the invertible operator
$S+\tau_0 B$, has a bounded inverse.
 By Lemma 2.3 in \cite{GMN99} and \eqref{x.5}, $M(z)$ is
invertible for
 all $z\in\bbC_+$, proving \eqref{x.6a}.

Since $S+zB$, $z\in \bbC_+, $ and $S+\tau_0 B$ are boundedly invertible
(dissipative) operators, Lemma~\ref{lx.0} implies
 $(\Xi(S+zB)-\Xi(S+\tau_0 B)) \in \calB_1(\calH)$.
By \eqref{x.5}, $(\Xi(S+\tau_0 B)-\Xi(S))\in
\calB_1(\calH)$
implying \eqref{x.6b}.
\end{proof}
\begin{remark}\lb{rkr.1}
We note that condition \eqref{ess} is crucial in connection
with
\eqref{x.5}. Indeed, there exists
 $S=S^*\in \calB(\calH)$
with
$
0\in \esss(S)
$
and a self-adjoint rank-one operator  $B$,
such that
\begin{equation}\lb{notin}
\Xi(S+B)-\Xi(S)=E_{S+B}((-\infty,0))-
E_{S}((-\infty,0))\notin \calB_1(\calH)
\end{equation}
 as implied by  a result of Krein \cite{Kr53}.
\end{remark}
\begin{corollary}\lb{ct.1}
Assume the hypotheses of Lemma~\ref{lx.1}. Then
\begin{equation}   \lb{c.1}
\tr (\Xi(S+tB)-\Xi(S))=0 \text{ for } t>0
\text{ sufficiently small}
\end{equation}
and
\begin{equation}\lb{c.2}
\tr (\Xi(S-tB)-\Xi(S))=\dim (\ker(S)) \text{ for } t>0
\text{ sufficiently small.}
\end{equation}
\end{corollary}
\begin{proof}
Eq. \eqref{c.1} follows from \eqref{x.6}
and Theorem~\ref{ttr.3}\,(iii),(iv). Next one notes that
for $t>0$ sufficiently small, $E_{S-tB}(\{0\})=0$ by
Remark~\ref{rx.1}, and
hence
\begin{align}
&\Xi(S-tB)-\Xi(S)=
E_{S-t B}((-\infty, 0))- E_S((-\infty,0)) \no \\
&=E_S([0, \infty))-E_{S-t B}([0, \infty))=
E_S([0, \infty))-E_{S-tB}((0, \infty)) \no \\
&=E_S(\{0\})+ E_S((0,\infty))-E_{S-tB}((0, \infty)) \no \\
&=E_S(\{0\})+ E_{-S}((-\infty, 0))-E_{-S+t B}((-\infty, 0))
\no \\
&=E_S(\{0\})+\Xi(-S) -\Xi(-S+tB).
\end{align}
Since by \eqref{c.1} (replacing $S$ by $-S$),
\begin{equation}\lb{c.3}
\tr (\Xi(-S) -\Xi(-S+tB))=0 \text{ for sufficiently small }
t>0,
\end{equation}
one obtains \eqref{c.2}.
\end{proof}

 Lemma~\ref{lx.0} yields $(\Xi(T)-\Xi(S))\in\calB_1(\calH)$
for $S=S^*\in\calB(\calH)$ and $T$ a dissipative
operator with $(T-S)\in \calB_1(\calH)$.

In the case of  general compact perturbations with trace
class imaginary parts, Hypothesis~\ref{hx.1} is not
sufficient
for the difference $(\Xi(S+A+iB)-\Xi(S))$ to be of
trace class
(see Remark~\ref{rce} below).
The following result, however,  shows that the pair
$(\Xi(S+A+iB),\Xi(S))$ has a trindex.

\begin{lemma}\lb{ltr.6}
Assume Hypothesis~\ref{hx.1}.~Then
\begin{equation}\lb{x.20}
(\Xi( S+A+iB) -\Xi( S+A))\in \calB_1(\calH)
\end{equation}
and the pair $(\Xi( S+A), \Xi( S))$ is a Fredholm
pair of orthogonal projections. Thus,
$(\Xi( S+A+iB), \Xi( S))$ has a trindex and
\begin{align}
&\trind (\Xi( S+A+iB) ,\Xi( S)) \no \\
&=\tr ( \Xi( S+A+iB) -\Xi( S+A)))+
\ind (\Xi(S+A), \Xi(S)) . \lb{tr.23}
\end{align}
\end{lemma}
\begin{proof}
Since by hypothesis $(S+A+\tau_0 B)^{-1}\in\calB(\calH)$
for some $\tau_0 \in \bbR$, \eqref{x.6b} implies
\eqref{x.20}.
Due to the fact that $0\notin \esss (S)$
(cf.~Remark~\ref{rx.1}) one can represent the spectral
projections $ \Xi(S+A)=E_{S+A}((-\infty,0))$ and
$\Xi(S)=E_{S}((-\infty,0))$ by the Riesz integrals
\eqref{x.7}, \eqref{x.8} and arguing as in the proof of
Lemma~\ref{lx.1} then yields
\begin{equation}\lb{x.21}
(\Xi(S+A)-\Xi(S))\in \calB_\infty(\calH).
\end{equation}
Thus, by Lemma~\ref{l2.1}\,(ii), the pair
$(\Xi(S+A),\Xi(S))$ is a
 Fredholm pair of orthogonal projections, which together
with \eqref{x.20} proves \eqref{tr.23}.
\end{proof}

\begin{lemma}\lb{lcom.2}
Let  $S$ be
 a signature operator, that is,   $S=S^*=S^{-1}$ and
$A=A^*\in \calB_2(\calH)$. Then
\begin{equation}\lb{ce.2}
(\Xi(S+A)-\Xi(S))\in \calB_1(\calH)
\end{equation}
if and only if
\begin{equation}\lb{comm}
[A,S]=(AS-SA)\in \calB_1(\calH).
\end{equation}
\end{lemma}
\begin{proof} Define two orthogonal projections
 $P=S_+$ and $Q=S_-$ such that $P+Q=I_\calH$ and
$
S=P-Q,
$
where $S_{\pm}=(|S|\pm S)/2$ (taking into
account that $\ker(S)=\{0\}$).
Next, we note
\begin{align}
&\frac{1}{2\pi i}\oint_{\Gamma}d\zeta
(S-\zeta)^{-1}A(S-\zeta)^{-1}
\no \\
&
=PAP\frac{1}{2\pi i}\oint_{\Gamma}d\zeta(1-\zeta)^{-2}
+QAQ\frac{1}{2\pi i}\oint_{\Gamma}d\zeta(-1-\zeta)^{-2}\no \\
&+(PAQ+QAP)\frac{1}{2\pi i}\oint_{\Gamma}d\zeta(1-\zeta)^{-1}
(-1-\zeta)^{-1}
=\frac{1}{4}S[S,A],  \lb{ce.4}
\end{align}
where the clockwise oriented contour $\Gamma$
encircles $\text{spec}(S+A)\cap
(-\infty,0)$.
On the other hand,
\begin{align}
&\Xi(S+A)-\Xi(S)=E_{S+A}((-\infty,0))
-E_{S}((-\infty,0))\no \\
&=-\frac{1}{2\pi i}\oint_{\Gamma}d\zeta (S-\zeta)^{-1}
A(S-\zeta)^{-1}\no \\
&\quad +\frac{1}{2\pi i}\oint_{\Gamma}d\zeta
(S+A-\zeta)^{-1}A(S-\zeta)^{-1}A(S-\zeta)^{-1}\no \\
&=-\frac{1}{4}S[S,A]+\frac{1}{2\pi i}\oint_{\Gamma}d\zeta
(S+A-\zeta)^{-1}A(S-\zeta)^{-1}A(S-\zeta)^{-1}
\lb{ce.3}
\end{align}
using \eqref{ce.4}.
Since $A\in \calB_2(\calH)$, the last term in \eqref{ce.3}
is a trace class operator and, hence,
$(\Xi(S+A)-\Xi(S))\in \calB_1(\calH)$ if and only
if $[S,A]\in
\calB_1(\calH)$.
\end{proof}

Combining the results of Lemmas~\ref{ltr.6}
 and \ref{lcom.2}
 we get the following result.
\begin{corollary}\lb{iff}
Assume $A=A^*\in \calB_2(\calH)$, $0\le B\in
\calB_1(\calH)$, and
 $S=S^*=S^{-1}$.  Then
\begin{equation}\lb{ce.5}
(\Xi(S+A+iB)-\Xi(S))\in \calB_1(\calH)
\end{equation}
if and only if
\begin{equation}\lb{ce.6}
[A,S]\in \calB_1(\calH).
\end{equation}
\end{corollary}

\begin{remark}\lb{rce}
Corollary~\ref{iff} illustrates why Hypothesis~\ref{hx.1}
is insufficient to garantee  \eqref{ce.5}. Indeed,
 choosing $\calH=\calK \oplus \calK$,  $S=I_\calK
\oplus (-I_\calK)$, and $A=\left(\begin{smallmatrix}
0& a \\ a & 0
\end{smallmatrix}\right)$, where
$a\in\calB_2(\calK)\backslash
\calB_1(\calK)$  with $||a||<1/2$ for some (necessarily
infinite
dimensional)
complex separable Hilbert space $\calK$, one  infers the
tripple $(S,A,B)$
 satisfies
Hypothesis~\ref{hx.1} for any $0\le B\in\calB_1(\calH)$  and
 $[A,S]\notin\calB_1(\calH)$.
Consequently, in this case,
\begin{equation}\lb{ce.7}
(\Xi(S+A+iB)-\Xi(S))\notin \calB_1(\calH)
\end{equation}
by Corollary~\ref{iff}.
\end{remark}

The next result concerns the continuity of
$\trind(\Xi( S+A+iB), \Xi(S))$ under small perturbations
of $A$ and $B$ in the operator and trace norm topology,
respectively.
\begin{theorem}\lb{ttr.7}
Assume Hypothesis~\ref{hx.1}. Let
$A_j=A^*_j\in \calB_\infty(\calH)$,
$0\leq B_j \in \calB_1(\calH)$, $j\in \bbN$, such that
$\lim_{j\to\infty}\|A_j-A\|=0$ and
$\lim_{j\to\infty}\|B_j-B\|_{\calB_1(\calH)}=0$. Then
\begin{equation}
\lim_{j\to \infty}\trind (\Xi(S+A_j+iB_j),
 \Xi(S)))=\trind (\Xi(S+A+iB),
 \Xi(S))).\lb{tr.26}
\end{equation}
\end{theorem}
\begin{proof}
By hypothesis, $(S+A+\tau_0 B)^{-1}\in\calB(\calH)$
for some $\tau_0\in\bbR$. Since
$\|A_j- A\|\to 0$ and $\|B_j -B\|\to 0$ as $j\to \infty$,
the
operators $S+A_j+\tau_0 B_j$ are also invertible for
$j\ge j_0$, $j_0$  sufficiently large. Since
the operator-valued logarithm is a continuous function of
 its (dissipative) argument in the $\calB(\calH)$-topology,
\begin{equation}\lb{tr.27}
\nlim_{j\to \infty} \log(S+A_j+\tau_0 B_j)=
\log(S+A+\tau_0 B)
\end{equation}
and hence by \eqref{x.2},
\begin{equation}\lb{tr.28}
\nlim_{j\to \infty} \Xi(S+A_j+\tau_0 B_j)=
 \Xi(S+A+\tau_0 B).
\end{equation}
By Remark~\ref{rx.1}, $0\notin \esss (S)$ and
thus
\begin{equation}\lb{tr.29}
(\Xi(S+A_j+\tau_0 B_j)-\Xi(S))\in \calB_\infty(\calH)
\end{equation}
and
\begin{equation}\lb{tr.30}
(\Xi(S+A+\tau_0 B)-\Xi(S))\in \calB_\infty(\calH).
\end{equation}
Thus the difference $(\Xi(S+A_j+\tau_0 B_j)
-\Xi(S+A+\tau_0 B))$ is a compact operator and hence
by \eqref{tr.28} and Lemma~\ref{ltr.4} one gets
\begin{equation}
\lim_{j\to \infty}\ind (\Xi(S+A_j+\tau_0 B_j)
,\Xi(S))=\ind (\Xi(S+A+\tau_0 B)
,\Xi(S))). \lb{tr.31}
\end{equation}
The estimate ($t\geq 0$)
\begin{align}
& \|(S+A+iB+it)^{-1}-(S+A+\tau_0 B+it)^{-1} \no \\
&-(S+A_j+iB_j+it)^{-1}+(S+A_j+\tau_0 B_j+it)^{-1}
\|_{\calB_1(\calH)} \no \\
&=(1+t^2)^{-1}o(1) \text{ as } j\to \infty,  \lb{tr.33}
\end{align}
with remainder term $o(1)$ uniform with respect to
$t\geq 0$,
then yields
\begin{align}
\lim_{j\to \infty}\| &\log(S+A_j+iB_j)-\log(S+A_j
+\tau_0 B_j)
\no \\
& -\log(S+A+iB)-\log(S+A+\tau_0
B)\|_{\calB_1(\calH)}=0. \lb{tr.32}
\end{align}
Hence
\begin{align}
\lim_{j\to \infty}
\|&\Xi(S+A_j+iB_j)-\Xi(S+A_j+\tau_0 B_j)  \no \\
& -\Xi(S+A+iB)+\Xi(S+A+\tau_0
B)\|_{\calB_1(\calH)}=0. \lb{x.32}
\end{align}
Combining \eqref{tr.31} and \eqref{x.32} yields
\begin{align}
\lim_{j\to \infty} &(\tr (\Xi(S+A_j+iB_j)-
\Xi(S+A_j+\tau_0 B_j))
+\ind ( \Xi(S+A_j+\tau_0 B_j), \Xi(S))) \no \\
&= \tr (\Xi(S+A+iB)-\Xi(S+A+\tau_0 B))
+\ind ( \Xi(S+A+\tau_0 B), \Xi(S)) \no \\
&=\trind (\Xi(S+A+iB),\Xi(S)),
\end{align}
proving \eqref{tr.26}.
\end{proof}

\section{Averaged Fredholm Indices and the Birman--Krein Formula} \lb{s4}

Assume Hypothesis~\ref{hx.1} with $A=0$.
By Remark~\ref{rx.1} the  operators $S+zB$, $z\in \bbC$
have a bounded inverse except for $z$ in a discrete set
$ \calD\subset \bbR$ with $\pm\infty$ the only possible
accumulation points of $\calD$. Introducing the family
of normal trace class operators
\begin{equation}\lb{t.45}
\calA(z)=B^{1/2}(S+zB)^{-1}B^{1/2},
 \quad z \in \bbC \backslash \calD.
\end{equation}
one verifies the resolvent-type identities
\begin{equation}\lb{t.46}
\calA(z_1)-\calA(z_2)
=(z_2 -z_1) \calA(z_1)\calA(z_2), \quad
\f{d}{dz}\calA(z)=-\calA(z)^2.
\end{equation}
Thus, $\calA(z_1)$ and
$\calA(z_2)$ commute for all $z_1,z_2\in\bbC\backslash
\calD$
and have a common complete orthogonal  system of
eigenvectors,
denoted by $\{\varphi_k \}_{k=1}^\infty$.

Making use of \eqref{t.46} one immediately gets the
following
result.
\begin{lemma}\lb{lt.4}
The operator $\calA(z)$, $z\in\bbC\backslash \calD$ has the
eigenvalue $z^{-1}$ with multiplicity
$\dim (\ker(S))$. In particular, if $S^{-1}\in\calB(\calH)$,
then
\begin{equation}
z^{-1} \notin \spec(\calA(z)). \lb{4.3}
\end{equation}
Let $\varphi$ be an eigenvector of $\calA(z_1)$
corresponding to the eigenvalue $\mu(z_1)$.
Then $\varphi$ is an eigenvector of $\calA(z_2)$
corresponding to the eigenvalue
\begin{equation}\lb{t.47}
\mu(z_2)=\frac{\mu(z_1)}
{1-\mu(z_1) (z_1-z_2)}
\end{equation}
of the same multiplicity and,
\begin{equation}\lb{t.48}
1-\mu(z_1) (z_1-z_2)\ne 0.
\end{equation}
\end{lemma}

Under Hypothesis~\ref{hx.1} with $A=0$, the operator
$\calA(\varepsilon)$
is well-defined for small (real) values of $\varepsilon$,
$\varepsilon \ne 0$,
$\calA(\varepsilon)\in \calB_1(\calH)$
(cf.~Remark~\ref{rx.1}) and hence
$\arctan (\calA(\varepsilon))\in \calB_1(\calH)$.
\begin{theorem}\lb{tt.bb}
Assume Hypothesis~\ref{hx.1} with $A=0$. Then
\begin{equation}\lb{lt.6-}
\lim_{\varepsilon \downarrow 0}
\tr( \arctan (\calA(\varepsilon) ))
=\int_\bbR \frac{dt\,n^*(t)}{1+t^2},
\end{equation}
where
 $n^*(t)$ is the left-continuous function.
\begin{equation}\lb{t.66a}
n^*(t)=\begin{cases} \sum_{s\in [0,t)}\,\dim(\ker (sB-S)),
&t >0, \\
0, &t=0, \\
-\sum_{s\in [t,0)}\dim(\ker (sB-S)), &t < 0. \end{cases}
\end{equation}
\end{theorem}
\begin{proof}
Fix a $\delta_0\in \bbR$, $\delta_0\notin\calD$. Denote by
$\{\mu_k(\delta_0)\}_{k\in\bbN}\subset
\spec (\calA(\delta_0))\backslash\{\delta_0^{-1}\}$ the
eigenvalues of $\calA(\delta_0)$  different from
$\delta_0^{-1}$ with corresponding multiplicities
$\{ m_k\}_{k\in\bbN}$.
First we  prove the following representation,
\begin{equation}\lb{t.49}
\lim_{\varepsilon \downarrow 0}
\tr( \arctan (\calA(\varepsilon)))=
\sum_{k=1}^\infty  m_k\arctan
 (\lambda_k)+
\frac{\pi}{2}\dim (\ker (S)),
\end{equation}
with
\begin{equation}\lb{t.50a}
\lambda_k=\frac{\mu_k(\delta_0)}{1-\mu_k(\delta_0)\delta_0},
\quad k\in \bbN.
\end{equation}
Since
$\calA(\delta_0)\in\calB_\infty(\calH)$, there exists a
$\gamma >0$ such that
\begin{equation}
\spec(\calA(\delta_0))\cap (\delta_0^{-1} - \gamma,
\delta_0^{-1} + \gamma)
=\begin{cases}
\emptyset,  & \dim(\ker (S))=0, \\
 \{\delta_0^{-1} \}, & \dim(\ker (S))>0.
\end{cases} \lb{t.50}
\end{equation}
For sufficiently small values of $\varepsilon \in \bbR$,
$\varepsilon \ne 0$,
the self-adjoint operator $\calA(\varepsilon)$
is well-defined and  the trace of
 $\arctan(\calA(\varepsilon))$ reads
\begin{equation}\lb{t.51}
\tr( \arctan (\calA(\varepsilon)))=
\sum_{k=1}^\infty m_k \arctan(\mu_k(\varepsilon))+
\dim(\ker (S)) \arctan (\varepsilon^{-1}),
\end{equation}
where by Lemma~\ref{lt.4},
\begin{equation}\lb{t.52}
\mu_k(\varepsilon)=\frac{\mu_k(\delta_0)}{1-\mu_k(\delta_0)
(\delta_0-\varepsilon)}, \quad k\in\bbN.
\end{equation}
 For $\varepsilon<\frac{\gamma \delta_0}{2}
\| \calA(\delta_0) \|^{-1}$ one obtains
\begin{align}
&\bigg |\sum_{k=1}^\infty m_k
\arctan (\mu_k(\varepsilon))
- \sum_{k=1}^\infty m_k\arctan (\lambda_k ) \bigg | \no \\
&\leq \sum_{k=1}^\infty m_k |\mu_k(\delta_0) |\,
|(1-\mu_k(\delta_0)\delta_0)^{-1} -
(1-\mu_k(\delta_0)(\delta_0-\varepsilon))^{-1}| \no \\
&\leq 2 \varepsilon (\gamma\delta_0)^{-2}
\sum_{k=1}^\infty m_k
|\mu_k(\delta_0) |^2  \le 2 \varepsilon
(\gamma\delta_0)^{-2}
\|\calA(\delta_0)
\|^2_{\calB_2(\calH)} . \lb{t.54}
\end{align}
By \eqref{t.54}, using the dominated convergence theorem,
one can perform the limit
 $\varepsilon  \downarrow 0$ in \eqref{t.51} and upon
combining \eqref{t.50a} and \eqref{t.52} one arrives at
\eqref{t.49} using Lemma~\ref{lt.4}.

Since the multiplicities
$m_k=\dim (\ker (B^{1/2}(S+\delta_0 B)^{-1}B^{1/2}-
\mu_k(\delta_0)I_\calH))$
of the eigenvalues $\mu_k(\delta_0)$ can also be computed as
\begin{equation}\lb{t.61}
m_k=\dim (\ker (B-\lambda_k S)),
\end{equation}
where $\lambda_k$ are given by \eqref{t.50a}, the absolutely
convergent series $\sum_{k=1}^\infty m_k
\arctan(\lambda_k)$ can be represented as the
Lebesgue integral
\begin{equation}\lb{t.62}
\sum_{k=1}^\infty m_k\arctan(\lambda_k)=
\sum_{k=1}^\infty m_k\int_0^{\lambda_k}\frac{dt}{1+t^2}=
\int_\bbR \frac{dt\,m(t)}{1+t^2},
\end{equation}
where $m(t)$ is the following eigenvalue counting function
\begin{equation}
m(t)=\begin{cases}
\sum_{\lambda\in (t,\infty)}\,\dim(\ker
(B-\lambda S)),  &t> 0, \\
 0, & t=0, \\
-\sum_{\lambda\in (-\infty, t)}\dim(\ker
(B-\lambda S)),&t < 0. \end{cases}
 \lb{t.63}
\end{equation}
Making the change of variables $t\to 1/t$ (separately on
$(-\infty,0)$ and $(0,\infty)$)
\begin{equation}\lb{t.64}
\int_\bbR \frac{dt\, m(t)}{1+t^2}=\int_\bbR
\frac{dt\, m(1/t)}{1+t^2},
\end{equation}
one infers by  \eqref{t.62} and \eqref{t.49} that
\begin{equation}
\lim_{\varepsilon \downarrow 0}\tr (
\arctan  (\calA(\varepsilon)))
=\int_\bbR \frac{dt\, m(1/t)}{1+t^2}+
\dim (\ker(S))\int_0^\infty \frac{dt}{1+t^2}
=\int_\bbR \frac{dt\, n(t)}{1+t^2}, \lb{t.65}
\end{equation}
where
\begin{equation}\lb{t.66}
n(t)=m(1/t) + \dim(\ker(S))\frac{1+\sgn(t)}{2},
\quad t\ne 0.
\end{equation}
Combining \eqref{t.63}
and \eqref{t.66} one concludes that $n(t)=n^*(t)$ for a.e.
$t\in\bbR$, where  $n^*(t)$ is a left-continuous function
given by \eqref{t.66a}.
\end{proof}

The following result is one of the main technical tools
in our paper.

\begin{theorem}\lb{tt.2}
Let $S=S^*\in \calB(\calH)$, $S^{-1}\in \calB(\calH)$,
$0\leq B\in \calB_1(\calH)$ and introduce
\begin{equation}\lb{t.6}
T(z)=S+zB, \quad z\in\ol {\bbC_+}.
\end{equation}
Define $\calZ=\bbC_+\cup \{ x\in \bbR \, |
\, (S+\tau x B)^{-1} \in \calB(\calH) $
for all $\tau \in [0,1] \}$. Then
\begin{equation}\lb{t.7}
(\log(T(z))-\log(S)) \in \calB_1(\calH), \quad z\in \calZ
\end{equation}
and
\begin{equation}\lb{t.8}
\tr (\log(T(z))-\log(S))=
\sum_{k=1}^\infty m_k \int_0^1 d\tau\, \frac{z
\lambda_k}{1+\tau z \lambda_k},\quad z\in \calZ,
\end{equation}
where $\{\lambda_k \}_{k\in\bbN}\subset \bbR$ is the
set  of eigenvalues with associated multiplicities
$\{m_k\}_{k\in\bbN}$ of
 the self-adjoint trace class operator
 $B^{1/2}S^{-1}B^{1/2}$.
\end{theorem}
\begin{proof}
First one notes that $T(z)$, $z\in \calZ$ is invertible.
 For $z\in
\bbR$ this holds by hypothesis and
for $\Im(z) >0$ this holds since $S^{-1}\in\calB(\calH)$.
Thus, $\log(T(z))$ and $\log(S)$ are well-defined.
In the following let $z\in \calZ$.  Since
 $B\in \calB_1(\calH)$, the representation
\begin{equation}\lb{t.11}
\log(T(z))-\log(S)=
-i\int_0^\infty dt\,((S+z B+it)^{-1}-(S+it)^{-1}),
\end{equation}
 and estimates for the resolvents $(S+iB+it)^{-1}$ and
$(S+it)^{-1}$, $t\ge 0$, analogous to those in the
proof of Lemma 2.6 of \cite{GMN99} yield
\eqref{t.7}. Therefore,
\begin{equation}
\tr(\log(T(z))-\log(S))
=-i\tr \bigg(\int_0^\infty dt\,((S+ z B+it)^{-1}
-(S+it)^{-1})\bigg ).
 \lb{t.12}
\end{equation}
Based on the estimates in the proof of Lemma~2.6
in \cite{GMN99}
mentioned above and the fact that $T=\int_0^\infty ds \,
T(s)\in\calB_1(\calH)$ and
\begin{equation}
\tr (T)=\int_0^\infty ds\,\tr( T(s)), \lb{t.12a}
\end{equation}
whenever $T(s)$
is continuous with respect to $s\in [0,\infty)$ in
$\calB_1(\calH)$-topology and $\|T(s)\|_{\calB_1(\calH)}
\leq C(1+s)^{-1-\varepsilon}$ for some $\varepsilon >0$,
one can interchange the integral and the trace in
\eqref{t.12} and obtain
\begin{equation}
\tr(\log(T(z))-\log(S))
=-i\int_0^\infty dt \,\tr((S+z B+it)^{-1}
-(S+it)^{-1})).
 \lb{t.13}
\end{equation}
Next, using the fact that $((S+\tau zB+it)^{-1}
-(S+it)^{-1})$
is differentiable with respect to $\tau$ in trace norm for
$(\tau,t)\in [0,1]\times [0,\infty)$ and
\begin{align}
&(d/d\tau)((S+\tau zB+it)^{-1}-(S+it)^{-1}) \no \\
&=-(S+\tau zB+it)^{-1}zB(S+\tau zB+it)^{-1} \lb{4.29}
\end{align}
in trace norm, one concludes
\begin{align}
&(d/d\tau)\tr((S+\tau zB+it)^{-1}-(S+it)^{-1}) \no \\
&=-\tr((S+\tau zB+it)^{-1}zB(S+\tau zB+it)^{-1}) \no \\
&=-\tr((S+\tau zB+it)^{-2}zB) \lb{4.30}
\end{align}
and hence
\begin{equation}
\tr ((S+z B+it)^{-1}-(S+it)^{-1})=
-\int_0^1d\tau \tr(S+\tau
zB+it)^{-2}zB), \quad t\ge 0,
 \lb{t.14}
\end{equation}
integrating \eqref{4.30} from $0$ to $1$ with respect to
$\tau$. Combining
\eqref{t.13} and
\eqref{t.14} one obtains
\begin{equation}\lb{t.19}
\tr (\log(T(z))-\log(S))=i\int_0^\infty dt
\int_0^1 d\tau\,\tr((S+\tau zB+it)^{-2}zB).
\end{equation}
Using the estimate
\begin{equation}
|\tr((S+\tau zB+it)^{-2}B)| \le
 \|(S+\tau zB+it)^{-1}\|^2 \|B\|_{\calB_1(\calH)}
 \le C(1+t^2)^{-1}, \lb{t.20}
\end{equation}
which holds uniformly with respect to $\tau\in [0,1]$,
Fubini's  theorem implies
\begin{equation}\lb{t.21}
\tr (\log(T(z))-\log(S))=i\int_0^1 d\tau\,
\int_0^\infty dt\,\tr((S+\tau zB+it)^{-2}zB).
\end{equation}
Applying  \eqref{t.12a} again, \eqref{t.20} implies
\begin{equation}
\int_0^\infty dt
\tr((S+\tau z B+it)^{-2}z B)
=\tr\bigg(\int_0^\infty dt\,
(S+\tau z B+it)^{-2} z B\bigg).\lb{t.22}
\end{equation}
Using
\begin{equation}\lb{t.23}
\int_0^\infty dt\,
(S+\tau zB+it)^{-2}=-i(S+\tau zB)^{-1},
\end{equation}
and combining \eqref{t.21}--\eqref{t.23} one finally gets
\begin{align}
\tr (\log(T(z))-\log(S))&=\int_0^1 d \tau\,
\tr ((S+\tau zB)^{-1}z B)
\no \\
& =\int_0^1 d\tau\, \tr (z B^{1/2} (S+\tau
z B)^{-1} B^{1/2}). \lb{t.24}
\end{align}

The trace of $zB^{1/2}(S+\tau zB)^{-1}B^{1/2}$,
$\tau\in [0,1]$, can easily be computed in terms of the
eigenvalues
$\{ \lambda_k \}_{k\in\bbN}$ of
$B^{1/2}S^{-1}B^{1/2}$. By Lemma~\ref{lt.4},
 $zB^{1/2}(S+\tau zB)^{-1}B $ has the eigenvalues
\begin{equation}\lb{t.25}
\mu_k(\tau,z)=\frac{z\lambda_k}{1+\tau z\lambda_k},
\quad k\in \bbN,
\end{equation}
with associated  multiplicities $\{m_k\}_{k\in\bbN}$. By
Lidskii's theorem (cf.~\cite[Ch.~3]{Si79})
\begin{equation}\lb{t.26}
\tr(zB^{1/2} (S+\tau zB)^{-1} B^{1/2})
 =\sum_{k=1}^\infty m_k\mu_k(\tau,z).
\end{equation}
By \eqref{t.24} and \eqref{t.25}
\begin{equation}\lb{t.27}
\tr(\log(T(z))-\log(S))=
\int_0^1d \tau \sum_{k=1}^\infty m_k\mu_k(\tau,z).
\end{equation}
Since $B^{1/2}S^{-1}B^{1/2}\in\calB_1(\calH)$ is
self-adjoint, one concludes that $\{\lambda_k\}_{k\in
\bbN}\subset \bbR$ and that the series $\sum_{k=1}^\infty
|\lambda_k|$  converges. Hence, applying the dominated
convergence theorem, one can interchange the sum and the
integral in \eqref{t.27}, arriving at \eqref{t.8}.
\end{proof}

\begin{corollary}\lb{ctt.2}
Under the assumptions of Theorem~\ref{tt.2},
\begin{equation}
\tr( \Im (\log(S+iB)-\Im(\log(S)))=
\tr( \arctan (B^{1/2} S^{-1} B^{1/2})).
\end{equation}
\end{corollary}
\begin{proof}
Pick $z=i$ in Theorem~\ref{tt.2}.
Taking the imaginary part of both sides of \eqref{t.8},
an explicit computation of  the integrals in \eqref{t.8}
 yields
\begin{equation}
\Im( \tr(\log(S+iB))-\tr(\log(S)))=
\sum_{k=1}^\infty m_k \arctan (\lambda_k)=
\tr(\arctan (B^{1/2}S^{-1}B^{1/2})).
\end{equation}
\end{proof}

\begin{remark}
Corollary~\ref{ctt.2} is
 an operator analog of the following elementary
fact
\begin{equation}\lb{tttt.1}
\Im(\log (a+ib))-\Im (\log(a))=
\arctan(b^{1/2}a^{-1}b^{1/2}), \quad
a\in\bbR\backslash\{0\}, \, b>0,
\end{equation}
where $\log(\cdot)$ and $\arctan(\cdot)$ denote the
corresponding principal branches, that is,
\begin{equation}\lb{tttt.2}
\Im(\log(\lambda))=0, \quad
  \lambda>0  \text{ and } -\frac{\pi}{2}< \arctan
(\lambda) <\frac{\pi}{2}, \, \lambda \in \bbR.
\end{equation}
\end{remark}

Theorem~\ref{tt.2} for $z=1$,
 has important consequences when computing the
trace of $(\Xi(S+B)-\Xi(S))$ (the case of self-adjoint
perturbations).

We start with the simplest  case of  self-adjoint
perturbations where $(S+tB)^{-1}\in\calB(\calH)$ for all
$t\in [0,1]$.

\begin{lemma}\lb{lt.2a}
Let $S=S^*\in \calB(\calH)$,
$0\leq B\in \calB_1(\calH)$. Assume
\begin{equation}\lb{t.2a}
(S+tB)^{-1}\in \calB(\calH) \text{ for all }
t\in [0,1].
\end{equation}
Then
\begin{equation} \lb{t.4a}
\tr(\Xi(S+B)-\Xi(S))=0.
\end{equation}
\end{lemma}
\begin{proof}
Under  hypothesis \eqref{t.2a} one can apply
Theorem~\ref{tt.2}
for $z=1$. Thus, $\tr(\Xi(S+B)-\Xi(S))\in\bbR$
and  therefore
\eqref{t.4a} holds due to \eqref{x.2} and the fact
that the left-hand side of \eqref{t.8} is real.
\end{proof}

Next we relax the condition \eqref{t.2a} of
invertibility  of $S+tB$ for all
$t\in [0,1]$, still assuming, however, that $S+\tau_0 B$ has
a bounded inverse for some $\tau_0\in \bbR$. The
following result
is concerned with the situation where the map $t\mapsto
(S+tB)^{-1}$, $t\in [-1,1]$ is singular at some points.

\begin{theorem}\lb{tt.4c}
Let $S=S^*\in \calB(\calH)$, $0\leq B\in \calB_1(\calH)$
and assume $(S+\tau_0 B)^{-1}\in \calB(\calH)$ for some
 $\tau_0\in \bbR$. Then
\begin{align}
\tr(\Xi(S+B)-\Xi(S))&=
-\sum_{s\in (0,1]} \dim(\ker (S+sB)), \lb{t.2c} \\
\tr(\Xi(S-B)-\Xi(S))&=
\sum_{s\in (-1,0]} \dim(\ker (S+sB)). \lb{t.2cc}
\end{align}
\end{theorem}
\begin{proof}
Since by hypothesis, $(S+\tau_0 B)^{-1}\in\calB(\calH)$,
Remark~\ref{rx.1} implies that
$(S+t B)^{-1}\in \calB(\calH)$ for all $t\in  [0,1]$
except possibly at a finite number
of points $0=t_0 < t_1< t_2 < ... <t_n<t_{n+1}=1$.

Introducing the notation
\begin{equation}\lb{t.3c}
E(t)=\Xi(S+tB), \quad t\in[0,1]
\end{equation}
one obtains for $\delta>0$ sufficiently small,
\begin{align}
&E_{S+B}((-\infty, 0))-E_S((-\infty, 0))=E(1)-E(0) \no \\
&=\sum_{k=1}^{n+1} (E(t_k)-E(t_k-\delta))+
\sum_{k=1}^{n+1}(E(t_k-\delta)-E(t_{k-1}+\delta)) \no \\
&+\sum_{k=1}^{n+1}(E(t_{k-1}+\delta)-E(t_{k-1})).
\lb{t.4c}
\end{align}
By Lemma~\ref{lt.2a},
\begin{equation}\lb{t.5c}
\tr \bigg(\sum_{k=1}^{n+1}(E(t_k-\delta)-E(t_{k-1}
+\delta))\bigg)=0,
\end{equation}
and by Corollary~\ref{ct.1} (for $\delta >0$ sufficiently
small),
\begin{equation}\lb{t.6c}
\tr\bigg(\sum_{k=1}^{n+1}(E(t_{k-1}+\delta)
-E(t_{k-1}))\bigg)=0,
\end{equation}
while
\begin{align}
&\tr \bigg(\sum_{k=1}^{n+1} (E(t_k)-E(t_k-\delta))\bigg)=
-\sum_{k=1}^{n+1}\dim(\ker(S+t_kB)) \no \\
&=-\sum_{s\in (0,1]}\dim(\ker (S+sB)).\lb{t.7c}
\end{align}
Combining \eqref{t.4c}--\eqref{t.7c} proves \eqref{t.2c}.

Setting $W=S-B$ one gets by \eqref{t.2c},
\begin{align}
&\tr\big( \Xi(S)-\Xi(S-B) \big )=
\tr\big(\Xi(W+B)-\Xi(W) ) \no \\
&=-\sum_{s\in (0,1]} \dim(\ker (W+sB))=
-\sum_{s\in (0,1]} \dim(\ker (S+(s-1)B)) \no \\
&=-\sum_{s\in (-1,0]}\dim(\ker (S+sB)),
\lb{t.8c}
\end{align}
proving \eqref{t.2cc}.
\end{proof}

As an immediate consequence one has the following result.

\begin{corollary}\lb{ctt.1}
Assume the hypotheses of Theorem~\ref{tt.4c}. Then
\begin{align}
&\tr(\Xi(S+tB)-\Xi(S)) \no \\
&=\ind (\Xi(S+tB),\Xi(S))=
\begin{cases}
 -\sum_{s\in (0,t]}\dim(\ker (S+sB)), &t> 0, \\
 0,  &t=0, \\
 \sum_{s\in (t,0]}\,\dim(\ker (S+sB))   ,&t < 0.
\end{cases}\lb{t.9c}
\end{align}
\end{corollary}
\begin{remark}\lb{safr}
The trace formula \eqref{t.9c}
 can be interpreted as follows. The Fredholm index of the pair
of spectral projections $(\Xi(S),\Xi(S+B))$
 coincides with the number of eigenvalues of
 $S+sB$ which cross the point
$0$ from the left  to the right  as the coupling
constant $s$ increases from $0$ to $1$.
\end{remark}

Now, we are prepared to prove our principal result.

\begin{theorem}\lb{ttr.8}
Assume Hypothesis~\ref{hx.1}.
Then the pair $( \Xi( S+A+iB),
\Xi(S))$
has a trindex and
\begin{equation}\lb{tr.35}
\trind
(\Xi( S+A+iB),
\Xi(S))=
\frac{1}{\pi}\int_\bbR \frac{dt\,n(t)}{1+t^2},
\end{equation}
where
\begin{equation}
n(t)=\ind  (\Xi(S+A+tB),\Xi(S)).
 \lb{tr.36}
\end{equation}
\end{theorem}
\begin{proof}
By Hypothesis~\ref{hx.1}, $(S+A+t B)^{-1}\in\calB(\calH)$
for all $t > 0$ sufficiently small. By
Corollary~\ref{ctt.2},
\begin{align}
&\tr( \Xi(S+A+tB +iB )-\Xi(S+A+t B)) \no \\
&=(1/\pi)\tr (\arctan
 ( B^{1/2}(S+A+t B)^{-1}B^{1/2})). \lb{y.1}
\end{align}
By Lemma~\ref{lx.1} one concludes that
\begin{equation}\lb{y.2}
\lim_{t \downarrow 0}\|\Xi(S+A+t B)
-\Xi(S+A) \|_{\calB_1(\calH)}=0
\end{equation}
and from standard properties of the operator logarithm
(cf.~\cite{GMN99}) one also infers
\begin{equation}\lb{y.3}
\lim_{t \downarrow 0}\|\Xi(S+A+iB+t B)
-\Xi(S+A+iB) \|_{\calB_1(\calH)}=0.
\end{equation}
Combining \eqref{y.1}--\eqref{y.3} one obtains
\begin{equation}
\tr(\Xi(S+A+iB)-\Xi(S+A))=(1/\pi)\lim_{t \downarrow 0}
 \tr( \arctan (B^{1/2}(S+A+t B)B^{1/2})) \lb{y.4}
\end{equation}
and by Theorem~\ref{tt.bb} one infers
\begin{equation}
\tr(\Xi(S+A+iB)-\Xi(S+A))=\frac{1}{\pi}
\int_\bbR \frac{dt\,n^*(t)}{1+t^2}=
\frac{1}{\pi}\int_\bbR \frac{dt\,n^*(-t)}{1+t^2}, \lb{y.5}
\end{equation}
where
\begin{equation}\lb{y.6}
n^*(-t)=\begin{cases}
\sum_{s\in [0,t)}\,\dim(\ker (sB-S-A)),  &t> 0, \\
0, &t=0, \\
-\sum_{s\in [t,0)}\dim(\ker (sB-S-A)),&t < 0. \end{cases}
\end{equation}
By Corollary~\ref{ctt.1},
\begin{equation}
n^*(t)=\ind(\Xi(S+A+tB),(\Xi(S+A)) \lb{4.67}
\end{equation}
and therefore,
\begin{equation}
\tr(\Xi(S+A+iB)-\Xi(S+A))=
\frac{1}{\pi}\int_\bbR dt\, \frac{\ind (\Xi(S+A+tB),
\Xi(S+A))}{1+t^2}.
\lb{y.7}
\end{equation}
Since by \eqref{x.21}, $(\Xi(S+A)-\Xi(A))\in
\calB_\infty(\calH)$, one concludes that
$(\Xi(S+A),\Xi(A))$ is a Fredholm pair of orthogonal
projections. Moreover, $(\Xi(S+A+iB)-\Xi(S))\in
\calB_1(\calH)$ implies that the pair
$(\Xi(S+A+iB),\Xi(S))$ has a trindex and hence
\begin{align}
&\trind (\Xi(S+A+iB),\Xi(S)) \no \\
&=\tr (\Xi(S+A+iB)-\Xi(S+A))+
\ind (\Xi(S+A), \Xi(S)). \lb{y.8}
\end{align}
Now \eqref{tr.35} follows from the chain rule
 \eqref{tr.7} for the index of a pair of projections,
\begin{align}
&\ind (\Xi(S+A+tB),\Xi(S+A))+\ind(\Xi(S+A),
\Xi(S)) \no \\
&=\ind (\Xi(S+A+tB),\Xi(S)),
\end{align}
and from the fact that the measure $(1/\pi)
(1+t^2)^{-1}dt$ is a probability measure on $\bbR$.
\end{proof}

\begin{remark} \lb{r4.10}
(i) If $B=0$, $n(t)$ is independent of $t$
and \eqref{tr.35} together with
$(1/\pi)\smallint_\bbR dt\,(1+t)^{-2}=1$ then imply
\begin{equation}
\trind(\Xi(S+A),\Xi(S))=\ind(\Xi(S+A),\Xi(S)), \lb{4.70}
\end{equation}
consistent with \eqref{2.19}.\\
(ii) The integral \eqref{tr.35} carries
out a ``smoothing'' of the integer-valued function $n(t)$
resulting in the expression of the trindex of a pair of
$\Xi$-operators. \\
(iii) Theorem~\ref{ttr.8} shows, in particular,
that under Hypothesis~\ref{hx.1}  the difference
\begin{equation}\lb{tr.37}
(\Xi( S+A+iB)-\Xi(S))\in \calB_1(\calH)
\end{equation}
if and only if
\begin{equation}\lb{tr.38}
(E_{S+A}((-\infty,0)) -
E_{S}((-\infty,0))) \in \calB_1(\calH).
\end{equation}
Under hypothesis \eqref{tr.38} one then obtains
\begin{align}
&\trind ( \Xi( S+A+iB),
 \Xi(S))=
\tr ( \Xi( S+A+iB)- \Xi(S)).
\lb{tr.39}
\end{align}
\end{remark}

\begin{remark}
Theorem~\ref{ttr.8} is an operator analog of the fact
\begin{align}
&\Im(\log (a+ib))-\Im (\log(a))=
\int_{\bbR}dt(1+t^2)^{-1}
(\chi_{(-\infty,0)}(a+tb)-\chi_{(-\infty,0)}(a)), \no \\
&\hspace*{8.5cm} a\in\bbR\backslash\{0\}, \, b>0,
\end{align}
where
$\chi_\Omega(\cdot)$ denotes the characteristic function
of $\Omega\subset\bbR$.
\end{remark}

There are two important special cases when
\eqref{tr.38} holds. For instance, if $S=I_\calH$
or $S=-I_\calH$ and $A=A^*\in\calB_\infty(\calH)$, the
difference \eqref{tr.38} is even a  finite-rank operator,
\begin{align}
\Xi(S+A)-\Xi(S) &=E_{S+A}((-\infty,0))
-E_{S}((-\infty,0)) \no \\
&=\begin{cases}
E_A((-\infty, -1)),  & S=I_\calH, \\
-E_A([1, \infty)), &  S=-I_\calH.  \end{cases}
\lb{tr.40}
\end{align}

\begin{lemma}\lb{ltrr.9}
Assume Hypothesis~\ref{hx.1} and $S=I_\calH$ or
$S=-I_\calH$. Then for $S=I_\calH$
\begin{equation}\lb{tr.42}
\Xi( I_\calH+A+iB)\in \calB_1(\calH)
\end{equation}
and
\begin{equation}\lb{tr.43}
\tr (\Xi( I_\calH +A+iB))=
\frac{1}{\pi}\int_\bbR \frac{dt\,n_-(t)}{1+t^2},
\end{equation}
where
\begin{equation}\lb{tr.44}
n_-(t)=\rank (E_{A+tB}((-\infty, -1)))
\end{equation}
is a decreasing right-continuous function.
 For $S=-I_\calH$ one has
\begin{equation}\lb{tr.45}
(\Xi( -I_\calH+A+iB)-I_\calH) \in \calB_1(\calH)
\end{equation}
and
\begin{equation}\lb{tr.46}
\tr(\Xi( -I_\calH +A+iB)-I_\calH)
=
-\frac{1}{\pi}\int_\bbR \frac{dt\,n_+(t)}{1+t^2},
\end{equation}
where
\begin{equation}\lb{tr.47}
n_+(t)=\rank (E_{A+tB}([1,\infty)))
\end{equation}
is an  increasing right-continuous function.
\end{lemma}
\begin{proof}
Let $S=I_\calH$. Then $\Xi(S)=\Xi(I_\calH)=0$ and
\begin{equation}
\ind ( \Xi(I_\calH +A+tB),\Xi(I_\calH) )=
\tr(\Xi(I_\calH +A+tB))
=\rank (E_{A+tB}((-\infty, -1))), \lb{tr.48}
\end{equation}
prove \eqref{tr.43} and \eqref{tr.44}. Next, let
$S=-I_\calH$. Then $\Xi(S)=\Xi(-I_\calH)=I_\calH$ and
\begin{align}
&\ind( \Xi(-I_\calH +A+tB), \Xi(-I_\calH))=
\tr(( \Xi(-I_\calH+A+tB)-I_\calH) \no \\
&=-\tr( E_{-I_\calH +A+tB}([0, \infty)))=
-\rank (E_{A+tB}([1,\infty))) \lb{tr.49}
\end{align}
prove \eqref{tr.46} and \eqref{tr.47}.
\end{proof}

Theorem~\ref{ttr.8} admits the following immediate
extension.

\begin{corollary}\lb{ctr.8}
Assume that the triples $(S,A_1,B_1)$ and $(S,A_2,B_2)$
satisfy
 Hypothesis~\ref{hx.1}. Then the pair $( \Xi( S+A_1+iB_1),
\Xi(S+A_2+iB_2))$ has a generalized trace and
\begin{equation}\lb{spr.35}
\gtr(\Xi( S+A_1+iB_1),\Xi((S+A_2+iB_2))=
\frac{1}{\pi}\int_\bbR \frac{dt\,n(t)}{1+t^2},
\end{equation}
where
\begin{equation}
n(t)=\ind  (\Xi(S+A_1+tB_1),\Xi((S+A_2+tB_2)).
 \lb{spr.36}
\end{equation}
\end{corollary}
\begin{proof}
By Theorem~\ref{ttr.8} the pairs $(\Xi( S+A_1+iB_1),
\Xi(S))$
and
$(\Xi( S+A_2+iB_2), \Xi(S))$ have a trindex and hence
the pair
$( \Xi( S+A_1+iB_1),\Xi(S+A_2+iB_2))$ has a generalized
trace and
\begin{align}
&\gtr(\Xi( S+A_1+iB_1),\Xi(S+A_2+iB_2)) \no \\
&=\trind( \Xi( S+A_1+B_1),\Xi(S))-\trind(\Xi( S+A_2+iB_2),
\Xi(S)). \lb{sprr.34}
\end{align}
Moreover, the following representations hold
\begin{align}
\trind (\Xi( S+A_1+iB_1),\Xi(S))&=
\frac{1}{\pi}\int_\bbR \frac{dt\,\ind  (\Xi(S+A_1+tB_1),
\Xi(S))}{1+t^2}, \lb{sprr.35} \\
\trind (\Xi( S+A_2+iB_2),\Xi(S))&=
\frac{1}{\pi}\int_\bbR \frac{dt\,\ind
(\Xi(S+A_2+tB_2),\Xi(S))}{1+t^2}. \lb{sprr.36}
\end{align}
Under Hypothesis~\ref{hx.1},
$(\Xi(S+A_j+tB_j)-\Xi(S))\in\calB_\infty(\calH)$, $j=1,2$
and
hence, by Theorem~\ref{ttr.3}\,(i),
\begin{align}
&\ind  (\Xi(S+A_1+tB_1),\Xi(S))-\ind
(\Xi(S+A_2+tB_2),\Xi(S))
\no \\
&=\ind  (\Xi(S+A_1+tB_1),\Xi(S+A_2+iB_2)). \lb{sprr.37}
\end{align}
Combining \eqref{sprr.34}--\eqref{sprr.37} proves
\eqref{spr.35} and \eqref{spr.36}.
\end{proof}

Finally, we turn to a  version of the Birman-Krein
formula \cite{BK62}.

\begin{theorem} \lb{ttr.10}
Under  Hypothesis~\ref{hx.1}, the operator
\begin{equation} \lb{tr.53}
{\bf S}=I_\calH -2iB^{1/2} (S+A+iB)^{-1}B^{1/2}
\end{equation}
is  unitary, $({\bf S}-I_\calH) \in \calB_1(\calH)$ and its
Fredholm determinant can be represented as follows
\begin{equation} \lb{tr.54}
\det ({\bf S})=\exp(-2\pi i\,\trind (\Xi(S+A+iB),\Xi(S))).
\end{equation}
\end{theorem}
\begin{proof}
Introduce the family of trace class operators
\begin{equation}\lb{tr.55}
\calA(z)=B^{1/2} (S+A+zB)^{-1}B^{1/2}, \quad \Im(z)>0.
\end{equation}
Then,
\begin{equation}\lb{tr.56}
{\bf S}=I_\calH-2i\calA(i), \quad ({\bf S}-I_\calH)\in
\calB_1(\calH)
\end{equation}
and hence the Fredholm determinant of ${\bf S}$ is
well-defined.
By the analytic Fredholm theorem the set of
 $z\in \bbC$ such that
$S+A+zB$ does not have a bounded inverse is discrete
and therefore there
exists a $\delta>0$ such that $\calA(\delta)$ is
well-defined.
By Lemma~\ref{lt.4},
the operator $\calA(\delta)$ has the eigenvalue
$\delta^{-1}$ if and only if
$\dim (\ker(S+A))\ne 0$ with associated multiplicity equal
to $\dim (\ker(S+A))$.

Let $\{\mu_k(\delta)\}_{k\in\bbN}=\spec(\calA(\delta))
\backslash\{\delta^{-1}\}$ different from the eigenvalue
$\delta^{-1}$ with multiplicities $\{m_k\}_{k\in\bbN}$.
By Lemma~\ref{lt.4}, $\calA(i)$ has the eigenvalues
\begin{equation}\lb{tr.59}
\mu_k(i)=\frac{\mu_k(\delta)}
{1-\mu_k(\delta) (\delta-i)}, \quad k\in\bbN,
\end{equation}
with multiplicities $\{m_k\}_{k\in\bbN}$ and, in addition,
the eigenvalue $-i$ of multiplicity
$\dim (\ker (S+A))$ (if $(S+A)$ has a nontrivial kernel).
Moreover, $\calA(i)$ has no other eigenvalues different
from zero.
Therefore, by \eqref{tr.56},
\begin{equation}\lb{tr.60}
\det ({\bf S})=(-1)^{\dim (\ker (S+A))}\prod_{k=1}^\infty
\bigg (\frac{1-\mu_k(\delta) \delta -i \mu_k(\delta)}
{1-\mu_k(\delta) \delta +i \mu_k(\delta)}
\bigg )^{m_k}.
\end{equation}
Moreover,
\begin{align}
&\trind (\Xi( S+A+iB), \Xi(S)) \no \\
&=\tr(\Xi( S+A+iB)- \Xi(S+A) )+\ind(\Xi( S+A),
 \Xi(S)) \no \\
&=(1/\pi) \lim_{\varepsilon \downarrow 0}
\tr(\arctan (B^{1/2} (S+A+\varepsilon B)^{-1} B^{1/2}))+
\ind(\Xi( S+A), \Xi(S)) \no \\
&= (1/\pi)\sum_{k=1}^\infty  m_k\arctan
(\mu_k(\delta)(1-\mu_k(\delta)\delta)^{-1})+
(1/2)\dim (\ker (S+A)) +n, \lb{tr.61}
\end{align}
for some $n\in\bbZ$. Combining \eqref{tr.60} and
\eqref{tr.61}
results in \eqref{tr.54}.
\end{proof}

To avoid additional technicalities we only treated
the case where $S$ is bounded. It is clear
that our formalism in Section~\ref{s3} extends to unbounded
dissipative operators $T$ in $\calH$, but such an
extension will be discussed elsewhere.

\section{Some Applications} \lb{s5}

The main purpose of this section is to obtain
new representations for Krein's spectral shift function
associated with a pair of self-adjoint operators
$(H_0,H)$  and to provide
 a generalization
of the classical Birman-Schwinger principle, replacing
the traditional eigenvalue counting functions by
appropriate spectral shift functions.

We start with our representation of Krein's spectral shift
function and temporarily assume the
following hypothesis.
\begin{hypothesis}\lb{h3.1}
Let  $H_0$ be a self-adjoint operator in $\calH$
with domain $\dom (H_0)$, $J$ a bounded self-adjoint
operator
with $J^2=I_{\calH}$, and $K\in \calB_2(\calH)$ a
Hilbert-Schmidt operator.
\end{hypothesis}

Introducing
\begin{equation}\lb{3.1}
V=KJK^*
\end{equation}
we define the self-adjoint operator $H$ in $\calH$ by
\begin{equation}\lb{3.2}
H=H_0+V, \quad \dom(H)=\dom(H_0).
\end{equation}
We could have easily incorporated the case where $K$ maps
between  different Hilbert spaces but we omit the
corresponding
details. Moreover, we introduce the following bounded
operator
in $\calH$,
\begin{equation}\lb{4.28}
\Phi (z)
=J+K^*(H_0-z)^{-1}K:\calH \rightarrow \calH, \quad
z\in \bbC\backslash \bbR.
\end{equation}
One easily verifies (cf.~\cite[Sect.~3]{GMN99}) that
$\Phi(z)$
is an operator-valued Herglotz function in $\calH$ (i.e.,
$\Im(\Phi(z))\geq 0$
for all $z\in\bbC_+$) and that
\begin{equation}\lb{3.8}
\Phi (z)^{-1}=J-JK^*(H-z)^{-1}KJ, \quad z\in\bbC
\backslash \bbR.
\end{equation}
In the following it is convenient to choose
\begin{equation}
J=\sgn (V),
\end{equation}
where in the present context the sign function is
defined by
$\sgn (x)=1$ if $x\geq 0$ and $\sgn (x)=-1$ if $x<0$.
Next, let $\{P_n\}_{n\in\bbN}$  be a family of finite-rank
spectral projections of $V$ satisfying,
\begin{equation}\lb{4.5}
\slim_{n\to \infty}P_n=I_{\calH}.
\end{equation}
Introducing the finite-rank operators
\begin{equation}
K_n=KP_n, \quad
V_n=(K_n)J(K_n)^*, \quad  n\in \bbN \lb{4.4},
\end{equation}
one infers (cf.~e.g., \cite{Gr73})
\begin{equation}\lb{4.5a}
\lim_{n\to \infty} \|V_n-V\|_{\calB_1({\calH})}.
\end{equation}
Together with the operator-valued Herglotz function
$\Phi(z)$ given by
\eqref{4.28}, we introduce the family of operator-valued
Herglotz functions
\begin{equation}\lb{4.37}
\Phi_n (z)
=J+K_n^*(H_0-z)^{-1}K_n:\calH \rightarrow \calH,
\end{equation}
and its finite-rank restriction
\begin{equation}\lb{4.38}
\Psi_n (z)=P_n\Phi_n(z)P_n :P_n\calH \rightarrow P_n\calH.
\end{equation}
One computes as in \eqref{3.8},
\begin{equation} \lb{4.39}
\Phi_n^{-1} (z)=J-JK_n^*(H_n-z)^{-1}K_nJ,
\end{equation}
where
\begin{equation} \lb{4.40}
H_n=H_0+V_n, \quad \dom (H_n)=\dom (H_0).
\end{equation}
Consequently,
\begin{equation}\lb{4.41}
\Psi_n^{-1} (z)=\Phi_n^{-1} (z)|_{P_n\calH} \text{ in }
 P_n\calH.
\end{equation}
Since $\Psi_n(z)$, $z\in \bbC_+$ is invertible in
$P_n\calH$,
one infers from \cite{GT97} (see also, \cite{Ca76}) the
existence of a family of  operators $\{\Xi_n (\lambda)\}$
defined for (Lebesgue) a.e.~$\lambda\in \bbR$, satisfying
\begin{equation}\lb{4.42}
0\le \Xi_{n}(\lambda)\le I_n
\text{ for a.e. } \lambda\in \bbR,
\end{equation}
where  $I_n$ denotes the identity operator in $P_n\calH$,
and
\begin{align}
\log (\Psi_n(z))&=C_n +\int_\bbR d\lambda \, \Xi_n(\lambda)
((\lambda-z)^{-1}-\lambda (1+\lambda^2)^{-1}), \quad
z\in \bbC\backslash \bbR, \lb{4.43} \\
C_n &=C_n^*=\Re(\log \Psi_n(i)), \lb{4.44} \\
\Xi_n(\lambda)&=\pi^{-1}\lim_{\varepsilon \downarrow 0}
\Im(\log (\Psi_n(\lambda+i\varepsilon ))
\text{ for a.e. } \lambda \in \bbR. \lb{4.45}
\end{align}

Next, we briefly recall the notion of Krein's spectral
shift
function for a pair of self-adjoint operators in $\calH$
(cf.~e.g., \cite[Sect.~19.1]{BW83}, \cite{BK62},
\cite{BP98},
\cite{BY93}, \cite{BY93a}, \cite{Ka78}, \cite{KS98},
\cite{Kr62},
\cite{Kr83}, \cite{Kr89}, \cite{KJ81}, \cite{Pu97},
\cite{Pu98a}, \cite{Si75}, \cite[Ch.~8]{Ya92} and the
literature cited therein), a concept originally
introduced by
Lifshitz \cite{Li52}, \cite{Li56}. Assuming $\dom(H_0)
=\dom(H)$ and $(H-H_0)\in\calB_1(\calH)$ (this could be
considerably relaxed but suffices for our present purpose),
Krein's (real-valued) spectral shift function
$\xi(\lambda,H_0,H)$ is uniquely defined for~a.e.
$\lambda\in\bbR$ by
\begin{align}
&\xi(\cdot,H_0,H)\in L^1(\bbR;d\lambda), \no \\
& \tr((H-z)^{-1}-(H_0-z)^{-1})=-\int_\bbR d\lambda\,
\xi(\lambda,H_0,H)(\lambda -z)^{-2}, \quad
z\in\bbC\backslash\bbR. \lb{4.45aa}
\end{align}

\begin{lemma}\lb{l4.7}
Denote by $\xi(\lambda,H_0,H_n)$ the spectral shift
function
associated with the pair $(H_0,H_n)$. Then
\begin{equation}
\xi(\lambda,H_0,H_n)=
\tr_{P_n\calH} ( \Xi_n(\lambda))-N_n \text{ for a.e. }
 \lambda \in \bbR, \lb{4.45a}
\end{equation}
where
\begin{equation}
N_n=\# \{ \text{of strictly negative eigenvalues of } V_n,
\text{\,counting multiplicity\,}\}. \lb{4.45b}
\end{equation}
\end{lemma}
\begin{proof} Differentiating $\tr_{P_n\calH}
(\log (\Psi_n(z)))$
with respect to $z$ (c.f., \cite[Sect.~IV.1]{GK69})
\begin{equation}\lb{4.46}
(d/dz)\tr_{P_n\calH}(\log(\Psi_n(z))=
\tr_{P_n\calH}(\Psi_n^{-1}(z)\Psi^{\prime}_n(z)),
\end{equation}
one obtains by \eqref{4.37} and \eqref{4.41}
\begin{equation}
\tr_{P_n\calH}(\Psi_n^{-1}(z)\Psi^{\prime}_n(z))=
\tr_{\calH}(P_n\Phi_n^{-1}(z)P_nP_n\Phi^{\prime}_n(z)P_n)=
\tr_{\calH}(\Phi_n^{-1}(z)\Phi^{\prime}_n(z)),
\lb{4.47}
\end{equation}
since $(d/dz)\Phi_n(z)=(d/dz)P_n\Phi_n(z)P_n$.
However, $\tr_{\calH}(\Phi_n^{-1}(z)\Phi^{\prime}_n(z))$
can be computed explicitly,
\begin{align}
& \tr_{\calH}((\Phi_n^{-1}(z)\Phi^{\prime}_n(z)) \no \\
&=\tr_{\calH}(J-JK_n^*(H_n-z)^{-1}K_nJ)
K_n^*(H_0-z)^{-2}K_n)
\no \\
&=\tr_{\calH}(K_n(J-JK_n^*(H_n-z)^{-1}
K_nJ)K_n^*(H_0-z)^{-2})
\no \\
&=\tr_{\calH}(H_0-z)^{-1}K_n(J-JK_n^*(H_n-z)^{-1}
K_nJ)K_n^*(H_0-z)^{-1})
\no \\
&=\tr((H_0-z)^{-1}K_nJK_n^*(H_0-z)^{-1}) \no \\
&-\tr_{\calH}(
(H_0-z)^{-1}K_nJK_n^*(H_n-z)^{-1}K_nJK_n^*(H_0-z)^{-1})
\no \\
&=\tr_{\calH}((H_0-z)^{-1}V_n(H_0-z)^{-1}-
(H_0-z)^{-1}V_n(H_n-z)^{-1}V_n(H_0-z)^{-1}) \no \\
&=-\tr_{\calH}( (H_n-z)^{-1}-(H_0-z)^{-1})
 \quad z\in \bbC\backslash \bbR, \lb{4.49}
\end{align}
iterating the second resolvent identity.
Combining \eqref{4.46}--\eqref{4.49} one infers
\begin{equation}\lb{4.50}
(d/dz)\tr_{P_n\calH}(\log(\Psi_n(z))=
-\tr_{\calH}( (H_n-z)^{-1}-(H_0-z)^{-1}).
\end{equation}
Taking traces in \eqref{4.43} one gets
\begin{equation}\lb{4.52}
\tr_{P_n\calH}(\log (\Psi_n(z)))=
\tr_{P_n\calH}(C_n) +
\int_\bbR d\lambda \,\tr_{P_n\calH}( \Xi_n(\lambda))
((\lambda-z)^{-1}-\lambda (1+\lambda^2)^{-1}),
\end{equation}
and thus, differentiating \eqref{4.52} with respect to $z$,
\begin{equation}\lb{4.53}
(d/dz)\tr(\log (\Psi_n(z)))=
\int_\bbR d\lambda \,\tr( \Xi_n(\lambda))
(\lambda-z)^{-2}.
\end{equation}
Comparing \eqref{4.53} and \eqref{4.50} one arrives at
the trace formula
\begin{equation}
\int_\bbR d\lambda \,\tr( \Xi_n(\lambda))
(\lambda-z)^{-2}=-\tr( (H_n-z)^{-1}-(H_0-z)^{-1}),
\quad z\in\bbC\backslash \bbR, \lb{4.54}
\end{equation}
and hence up to an additive constant, $\tr( \Xi_n(\lambda))$
coincides  with the spectral shift function
$\xi(\lambda,H_0,H_n)$ associated with the pair $(H_0,H_n)$
for a.e.~$\lambda \in \bbR$.

Next we determine this constant. Introducing
$J_n=J\big|_{P_n\calH}$, $J_n^2=I_n$, one obtains
\begin{align}
\log(J_n)&=i\Im(\log(J_n)) \no \\
&=\pi^{-1}\int_\bbR d\lambda \,\Im(\log(J_n) )
((\lambda-z)^{-1}-\lambda (1+\lambda^2)^{-1}). \label{4.43a}
\end{align}
Moreover,
\begin{equation}\label{4.43b}
\pi^{-1}\tr_{P_n\calH}(\Im
(\log ((J_n)))=N_n,
\end{equation}
where $N_n$ denotes the number of strictly negative
eigenvalues
of $J_n$, counting multiplicity. Thus, $N_n$ coincides
with the
number of strictly negative eigenvalues of $V_n$. Combining
\eqref{4.43a}, \eqref{4.43b}, and using \eqref{4.43}
results in
\begin{align}
&\log (\Psi_n(z))-\log(J_n) \no \\
&=C_n +\int_\bbR d\lambda \, (\Xi_n(\lambda)-\pi^{-1}\Im
(\log ((J_n)))((\lambda-z)^{-1}-\lambda (1+\lambda^2)^{-1})
\lb{4.44b}
\end{align}
and hence in
\begin{equation}
\tr_{P_n\calH}(\Im(\log (\Psi_n(z))-\log(J_n)))
=\int_\bbR d\lambda \,(\tr( \Xi_n(\lambda))-N_n)
\frac{\Im(z)}{(\lambda-\Re(z))^{^2}+\Im(z)^2}. \lb{4.56}
\end{equation}
Since $\vert\vert (H_0-iy)^{-1}\vert\vert =O(\vert y
\vert^{-1}) $ as $ y \uparrow +\infty$, one concludes
\begin{equation}\lb{4.57}
y\vert\vert
\log (\Psi_n(iy))-\log(J_n) \vert\vert=O
(1) \text{ as } y \uparrow + \infty
\end{equation}
and hence that
$y\Im(\tr_{P_n\calH} (\log(\Psi_n(iy))-\log (J_n)))$
is bounded $\text{ as }  y \uparrow +\infty$.
In particular, \eqref{4.56} and \eqref{4.57} imply that
\begin{equation}\lb{4.58}
\xi_n(\lambda)=\tr_{P_n\calH} (\Xi_n(\lambda))-N_n
\end{equation}
is integrable,
\begin{equation}\lb{4.60}
\xi_n(\cdot) \in L^1(\bbR;d\lambda).
\end{equation}
Since
\begin{equation}
\int_\bbR d\lambda(\lambda-z)^{-2}=0 \text{ for all }
z\in\bbC, \, \Im(z)\ne 0,
\end{equation}
\eqref{4.58} and \eqref{4.54} yield the trace formula
\begin{align}
&\int_\bbR d\lambda \,\xi_n(\lambda)(\lambda-z)^{-2}=
\int_\bbR d\lambda \,\tr_{P_n\calH}( \Xi_n(\lambda))
(\lambda-z)^{-2} \no \\
&=-\tr_{\calH}( (H_n-z)^{-1}-(H_0-z)^{-1}),\quad z\in
\bbC\backslash \bbR,
\lb{4.62}
\end{align}
which together with \eqref{4.60} proves
\eqref{4.45a}, \eqref{4.45b}.
\end{proof}

\begin{theorem}\lb{t5.8}
Assume Hypothesis~\ref{h3.1} and fix a $p>1$. Moreover,
let
$V=KJK^*$, where $J=\sgn (V)$. Then the spectral shift
function
$\xi(\lambda,H_0,H)$ associated with the pair $(H_0,H)$,
$H=H_0+V$ admits the representation
\begin{equation}
\xi(\lambda,H_0,H)=\trind(\Xi(J+K^*(H_0-
\lambda-i0)^{-1}K), \Xi(J)) \text{ for a.e. }
\lambda \in\bbR, \lb{5.52}
\end{equation}
where $K^*(H_0-\lambda-i0)^{-1}K$ is defined as
\begin{equation}
\lim_{\varepsilon\downarrow 0}\|K^*(H_0-\lambda-i0)^{-1}K-
K^*(H_0-\lambda-i\varepsilon)^{-1}K\|_{\calB_p(\calH)}=0
\text{ for a.e. }  \lambda \in \bbR. \lb{4.1}
\end{equation}
\end{theorem}
\begin{proof}
First of all, one notes that  the boundary values
$K^*(H_0-\lambda-i0)^{-1}K$ and $K^*(H-\lambda-i0)^{-1}K$
exist
$\lambda$ a.e. in the topology $\calB_p(\calH)$ for
 every $p>1$ (but in general not for $p=1$,) \cite{Na89},
\cite{Na90} (see also \cite[Ch.~3]{BW83},
 \cite{BE67}, \cite{de62}). By
\eqref{3.8}, the operator
$J+K^*(H_0-\lambda-i0)^{-1}K$ has a bounded inverse
for a.e.
$\lambda\in\bbR$. Moreover (see, e.g.,
\cite[Sect.~I.3.4]{BW83}),
\begin{align}
&\lim_{\varepsilon \downarrow 0} \|\Im(K^*(H_0-
\lambda-i\varepsilon)^{-1}K)-\Im(K^*(H_0-
\lambda-i0)^{-1}K)\|_{\calB_1(\calH)}=0  \lb{5.1} \\
 & \hspace*{8.6cm} \text{ for a.e. } \lambda \in \bbR. \no
\end{align}
Thus, there exists a set $\Lambda\subset\bbR$ with
$|\bbR\backslash\Lambda|=0$ ($|\cdot|$ denoting Lebesgue
measure on $\bbR$) satisfying the following properties.
For any $\lambda\in \Lambda:$

\noindent (i) The boundary values
$\Phi(\lambda+i0)=\lim_{\varepsilon
\downarrow 0}\Phi(\lambda+i\varepsilon)$ exist in
$\calB_p(\calH)$-topology (cf. \eqref{4.28}).

\noindent (ii) The operator $\Phi(\lambda+i0)$ has a
bounded
inverse.

\noindent (iii)  $\Im(\Phi(\lambda +i\varepsilon))$
converges
to $\Im(\Phi(\lambda+i0))$ as $\varepsilon \downarrow 0$
in $\calB_1(\calH)$-topology.

For any $n\in \bbN$, $\lambda\in \Lambda$ introduce
the function
\begin{equation}\lb{4.6}
\xi_n(\lambda)=\trind( \Xi(P_n\Phi(\lambda+i0)P_n,
 \Xi(P_nJP_n)).
\end{equation}
Since $P_n$ commute with $J$ and the subspace $P_n \calH$
is invariant for
$J+P_nK^*(H_0-\lambda-i0)^{-1}KP_n$ one
concludes by \eqref{4.6} that (cf. \eqref{2.20})
\begin{align}
\xi_n(\lambda)&=\tr ( \Xi(J+P_n
K^*(H_0-\lambda-i0)^{-1}KP_n)-\Xi(J)) \no \\
&=\trind ( \Xi(J+P_n
K^*(H_0-\lambda-i0)^{-1}KP_n), \Xi(J)), \quad
\lambda \in \Lambda. \lb{4.7}
\end{align}
On the other hand, by Lemma~\ref{l4.7}, \eqref{4.37},
 and \eqref{4.38} one infers
that the function $\xi_n(\lambda)$ coincides with the
spectral shift function
$\xi(\lambda, H_0, H_n)$ associated with the pair
$(H_0, H_n)$
\begin{equation}\lb{4.79}
\xi_n(\lambda)=\xi(\lambda, H_0, H_n),\text{ a.e. }
\lambda \in \bbR,
\end{equation}
 where $H_n$ is given
by \eqref{4.40}.
By a result of Gr\"umm \cite{Gr73}, properties
(i)--(iii), and \eqref{4.5}, one obtains for
$\lambda\in \Lambda$,
\begin{align}
&\lim_{n\to \infty}
\|(\Re (P_nK^*(H_0-\lambda-i0)^{-1}KP_n)
-\Re ( K^*(H_0-\lambda-i0)^{-1}K)\|_{
\calB_p(\calH)}=0, \lb{4.8} \\
&\lim_{n\to \infty}
\|(\Im ( P_nK^*(H_0-\lambda-i0)^{-1}KP_n)-\Im (
 K^*(H_0-\lambda-i0)^{-1}K)\|_{\calB_1(\calH)}=0.
\lb{4.9}
\end{align}
Applying the approximation Theorem~\ref{ttr.7} then yields
\begin{equation}\lb{4.10}
\lim_{n\to \infty} \xi_n(\lambda)=
\trind( \Xi(\Phi(\lambda+i0), \Xi(J)), \quad
\lambda \in \Lambda \text{ pointwise}.
\end{equation}
Convergence \eqref{4.5a} of $V_n$ to $V$ in trace norm
implies the convergence of the corresponding spectral
shift
 functions $\xi_n(\lambda)$
to the spectral shift function $\xi(\lambda,H_0,H)$ in
$L^1(\bbR)$. This in turn implies the existence of a
subsequence
$\{ \xi_{n_k}(\lambda)\}_{k\in\bbN}$ converging pointwise
a.e. Together with \eqref{4.10} this proves \eqref{5.52}.
\end{proof}

\begin{remark}\lb{rindex}
Let $\Lambda=\{\lambda\in\bbR\,|\, \text{s.\,t. } A(\lambda)
\text{ and } B(\lambda)
\text{ exist, } A(\lambda)\in\calB_2(\calH),
B(\lambda)\in\calB_1(\calH), \text{ and }
(J+A(\lambda)+iB(\lambda))^{-1}\in\calB(\calH) \}$. Then,
as shown in Corollary~\ref{iff}, the condition
\begin{equation}
[A(\lambda),J]\notin \calB_1(\calH) \lb{5.49}
\end{equation}
is necessary and sufficient for the validity of
\begin{equation}
(\Xi(J+A(\lambda)+iB(\lambda))-\Xi(J))\notin
\calB_1(\calH).
\lb{5.53}
\end{equation}
Next, note that the following three conditions,
\begin{align}
(i) \, &\rank (E_V((-\infty,0))=\rank (E_V((0,
 \infty))=\infty, \\
(ii)\, & \lambda\in \esss (H_0), \\
(iii)\, & A(\lambda)\notin \calB_1(\calH),
\end{align}
are a consequence of condition \eqref{5.49}.
Thus, if at least one of the conditions (i)--(iii) is
violated,
the $\xi$-function (c.f., \eqref{5.52}) can be
represented in
the simple form
\begin{equation}
\xi(\lambda,H_0,H)=\tr(\Xi(J+K^*(H_0-
\lambda-i0)^{-1}K)- \Xi(J)),
 \lb{rind.1}
\end{equation}
and the concept of a trindex becomes redundant in this case.
On the other hand, there are of course examples (c.f.,
Remark~\ref{rce}), where
\begin{equation}
|\{\lambda \in \bbR \, |\, (\Xi(J+|V|^{1/2}(H_0-
\lambda-i0)^{-1}|V|^{1/2})- \Xi(J))\notin
\calB_1(\calH) \}|
>0,
\lb{5.53a}
\end{equation}
with $|\cdot |$ denoting Lebesgue measure on $\bbR$. A
concrete
example illustrating \eqref{5.53a} can be constructed as
follows.
Consider an infinite dimensional complex separable
Hilbert space
$\calK$, an operator $0\leq
k=k^*\in\calB_2(\calK)\backslash\calB_1(\calK)$, with
$\ker(k)=\{0\}$, and a self-adjoint operator $h_0$ in
$\calK$ such
that
\begin{equation}
a(\lambda)=\nlim_{\varepsilon\to 0}\Re(k(h_0-\lambda
-i\varepsilon)^{-1}k)\notin\calB_1(\calK) \text{ for~a.e. }
\lambda\in\bbR. \lb{5.54}
\end{equation}
Existence of such $\calK$, $k$, and $h_0$ can be
inferred from
\cite{Na89}. Next, define $\calH=\calK\oplus\,\calK$ and
introduce $H_0=\left(\begin{smallmatrix} h_0& 0 \\ 0 & 0
\end{smallmatrix}\right)$, $V=
i\left(\begin{smallmatrix} 0& k^2
 \\ -k^2 & 0 \end{smallmatrix}\right)$. Then $J=\sgn(V)=
i\left(\begin{smallmatrix} 0& I_\calK \\ -I_\calK & 0
\end{smallmatrix}\right)$, $|V|^{1/2}=
\left(\begin{smallmatrix}
k & 0 \\ 0 & k \end{smallmatrix}\right)$, and
\begin{align}
A(\lambda)=\nlim_{\varepsilon\to
0}
\Re(|V|^{1/2}(H_0-\lambda-i\varepsilon)^{-1}|V|^{1/2})
&=\begin{pmatrix}a(\lambda) & 0 \\ 0
& -(1/\lambda)k^2 \end{pmatrix} \lb{5.55} 
\end{align}
for~a.e.~$\lambda\in\bbR\backslash\{0\}$. Moreover, one computes
\begin{equation}
[A(\lambda),J]=i\begin{pmatrix} 0& a(\lambda)
+(1/\lambda)k^2 \\ a(\lambda)+(1/\lambda)k^2 & 0
\end{pmatrix} \notin \calB_1(\calH) \text{ for~a.e. }
\lambda\in\bbR, \lb{5.56}
\end{equation}
since $k^2\in\calB_1(\calK)$.

\end{remark}

As a consequence of Theorems~\ref{ttr.8} and \ref{t5.8}
one has the following representation for the spectral
shift function via the integral of the index of a
Fredholm pair of spectral projections.

\begin{theorem} \lb{main}
Assume Hypothesis~\ref{h3.1} and introduce $V=KJK^*$. In
addition, for a.e.~$\lambda\in\bbR$, let
$A(\lambda)+iB(\lambda)$  be  the normal boundary values
of the
operator-valued Herglotz function $K^*(H_0-z)^{-1}K$
on the real axis, that is,
\begin{equation}
A(\lambda)=\nlim_{\varepsilon \downarrow 0}
 \Re(K^*(H_0-\lambda-i\varepsilon)^{-1}K)
 \text{ for a.e. } \lambda\in\bbR
\end{equation}
and
\begin{equation}
B(\lambda)=\nlim_{\varepsilon \downarrow 0}
\Im(K^*(H_0-\lambda-i\varepsilon)^{-1}K)
 \text{ for a.e. } \lambda\in\bbR.
\end{equation}
Then
\begin{align}
&\xi(\lambda, H_0,H)=\frac{1}{\pi} \int_\bbR dt
\frac{\ind \big(E_{J+A(\lambda)+tB(\lambda)}
((-\infty,0)),E_J((-\infty,0))\big )}{1+t^2}
\lb{pus} \\
&\hspace{8.1cm} \text{for  a.e. } \lambda\in\bbR. \no
\end{align}
\end{theorem}

 In the particular case of
sign-definite perturbations, that is,
$J=I_\calH$ or $J=-I_\calH$, applying Lemma~\ref{ltrr.9}
yields the following result originally due to
Pushnitski~\cite{Pu97}, representing the spectral shift
function in terms of an integrated eigenvalue
counting function.
\begin{corollary}\lb{c5.6} \mbox{\rm
(Pushnitski~\cite{Pu97}.)}
Let $0\le V\in \calB_1(\calH)$ and $H_0=H_0^*$. Then the
spectral shift
function $\xi(\lambda, H_0, H_0\pm V)$ associated
 with the pair $(H_0, H_0\pm V)$ admits the representation
\begin{equation} \lb{puss}
\xi(\lambda, H_0,H_0\pm V)=\pm\frac{1}{\pi} \int_\bbR dt
\frac{\rank (E_{\mp (A(\lambda)+
tB(\lambda)}((1, \infty)))}
{1+t^2}.
\end{equation}
\end{corollary}

\begin{remark} \lb{r5.5}
(i) Strictly speaking, a direct  application of
Lemma~\ref{ltrr.9} in the case
of nonpositive perturbations would give the representation
\begin{equation}
\xi(\lambda, H_0,H_0- V)=-\frac{1}{\pi} \int_\bbR dt
\frac{\rank (E_{ A(\lambda)+tB(\lambda)}
([1, \infty)))}
{1+t^2},
\end{equation}
which, however, yields the same result as in
\eqref{puss}, since
\begin{equation}
\int_\bbR dt
\frac{\rank (E_{ A(\lambda)+tB(\lambda)}(\{1\}))}
{1+t^2}=0 \text{ for a.e. }\lambda\in \bbR.
\end{equation}
(ii) In the special case where $\lambda \in
\bbR\backslash\{\spec (H_0)\cup \spec (H)\}$,
\eqref{pus} turns into
\begin{align}
\xi(\lambda, H_0,H)&=\ind(E_{J+A(\lambda)}((-\infty, 0)),
E_J((-\infty, 0)))
\no \\
&=\tr (E_{J+A(\lambda)}((-\infty, 0))-
E_J((-\infty, 0)))
. \lb{5.52a}
\end{align}
In the particular cases of sign-definite perturbations
($J=\pm I_\calH$) one obtains
\begin{equation}
\xi(\lambda, H_0,H_0\pm V)=\pm\text{rank}
 (E_{\mp A(\lambda)}((1, \infty))), \lb{sobo}
\end{equation}
since $\text{rank} (E_{\mp A(\lambda)}(\{1\}))=0$ for
a.e.~$\lambda\in\bbR$.
 The result \eqref{sobo} is due to Sobolev \cite{So93}.
\end{remark}

\vspace*{3mm}

The trindex representation \eqref{5.52} for $\xi(\lambda,H_0,H)$
enables us to introduce a new generalized spectral
 shift function outside the trace-class perturbation scheme
under rather weak assumptions on $H-H_0$. First we recall the
following exponential representation for operator-valued Herglotz
functions partially proven in \cite{GMN99}.

\begin{theorem}\lb{t5.8a}
Suppose $M:\bbC_+ \rightarrow \calB(\calH)$ is an
operator-valued Herglotz function
and $M(z_0)^{-1} \in \calB(\calH)$ for some {\rm (}and hence
for all{\rm  \,)} $z_0\in\bbC_+$. Then
there exists a family of bounded self-adjoint
weakly {\rm (}Lebesgue{\rm \,)} measurable operators
 $\{\what\Xi(\lambda) \}_{\lambda\in \bbR}\subset \calB(\calH)$,
\begin{equation}\lb{5.64}
0\le \what\Xi(\lambda)\le I_\calH \text{ for a.e. }
\lambda\in \bbR
\end{equation}
such that
\begin{equation}\lb{5.65}
\log(M(z))=C+
\int_\bbR d \lambda \, \what\Xi(\lambda)
((\lambda-z)^{-1}-\lambda(1+\lambda^2)^{-1}),
 \quad z\in \bbC_+
\end{equation}
the integral taken in the weak sense, where $C=C^*=\Re(\log(M(i))) \in
\calB(\calH)$.
Moreover, suppose there exists a measurable subset $\Lambda\in\bbR$,
$|\Lambda|\neq 0$ such that
\begin{equation}
\nlim_{\varepsilon\to 0}M(\lambda+i\varepsilon)
=M(\lambda+i0)\in\calB(\calH) \text{ for a.e. }
\lambda\in\Lambda \lb{5.66}
\end{equation}
such that $\Im(M(\lambda+i0))\geq 0$ and
$M(\lambda+i0)^{-1}\in\calB(\calH)$ for a.e.~$\lambda\in\Lambda$. Then
\begin{equation}
\what \Xi(\lambda) =\Xi(M(\lambda+i0)) \text{ for a.e. }
\lambda\in\Lambda. \lb{5.67}
\end{equation}
\end{theorem}
\begin{proof}
Since \eqref{5.64} and \eqref{5.65} have been proven in \cite{GMN99},
we focus on \eqref{5.67}. Let $\{P_n\}_{n\in\bbN}$ be an increasing
family of orthogonal projections of rank $n$, that is, $\rank(P_n)=n$,
$P_n\calH\subset P_{n+1}\calH$, with
\begin{equation}
\slim_{n\to \infty} P_n=I_\calH. \lb{5.67a}
\end{equation}
Combining the norm continuity of the logarithm of bounded dissipative
operators as discussed in Section~2 of \cite{GMN99} with the
 exponential Herglotz 
representation for $P_n\log(M(z))P_n$ (i.e., the finite-dimensional analog 
of \eqref{5.65} as in \eqref{4.43}--\eqref{4.45}), one infers for each
$n\in\bbN$ the existence of a subset $\Lambda_n\subseteq\Lambda$,
$|\Lambda\backslash\Lambda_n|=0$ such that 
\begin{equation}
P_n\what\Xi(\lambda)P_n=P_n\Xi(M(\lambda+i0))P_n \text{ for all }
\lambda\in\Lambda_n. \lb{5.67b}
\end{equation}
Thus,
\begin{equation}
P_n\what\Xi(\lambda)P_n=P_n\Xi(M(\lambda+i0))P_n \text{ for all }
\lambda\in\bigcap_{m\in\bbN}\Lambda_m, \,\, n\in\bbN. \lb{5.67c}
\end{equation}
Since $|\Lambda\backslash\bigcap_{n\in\bbN}\Lambda_n|=0$ one concludes
\eqref{5.67}.
\end{proof}

Next, assuming $H_0$ and $V$ to be self-adjoint operators in $\calH$
with corresponding domains $\dom (H_0)$ and $\dom(V)$, such
that
\begin{equation}
\dom(|V|^{1/2})\supseteq \dom(|H_0|^{1/2}), \lb{5.68}
\end{equation}
and introducing the signature operator
\begin{equation}
J=\sgn(V)  \text{ with } J|_{\ker(V)}=I_\calH|_{\ker(V)},
\lb{5.69}
\end{equation}
we may define the bounded operator-valued Herglotz function
$\phi(z)\in\calB(\calH)$
\begin{equation}
\phi(z)=J+\ol{|V|^{1/2}(H_0-z)^{-1}|V|^{1/2}}\, , \quad z\in\bbC_+.
\lb{5.70}
\end{equation}
(here the bar denotes the operator closure in $\calH$.)
In addition, we suppose that
\begin{equation}
\phi(z)^{-1}\in\calB(\calH) \text{ for some (and hence for all) }
z\in\bbC_+. \lb{5.70a}
\end{equation}
Applying Theorem~\ref{t5.8a}, one
concludes that $\phi(z)$ admits the representation
\begin{align}
&\log (\phi(z))=C +\int_\bbR d\lambda \, \what \Xi(\lambda,H_0,H)
((\lambda-z)^{-1}-\lambda (1+\lambda^2)^{-1}), \quad
z\in \bbC\backslash \bbR, \lb{5.71} \\
&C =C^*=\Re(\log \phi(i)), \quad
0\le \what\Xi(\lambda,H_0,H)\le I_\calH
\text{ for a.e. } \lambda\in \bbR. \lb{5.72}
\end{align}
Here our notation $\what\Xi(\lambda,H_0,H)$ emphasizes the underlying
pair $(H_0,H)$, where $H$``$=$''$H_0+V$ formally represents the
perturbation  of $H_0$ by $V$. We will return to a discussion of this
point in  Remark~\ref{r5.9a} below.

\begin{definition}\lb{d5.9}
In addition to \eqref{5.68}--\eqref{5.72} assume the existence of
a (Lebesgue) measurable set $\Lambda$, $|\Lambda|\neq 0$, such that the
pair
$(\what\Xi(\lambda,H_0,H),\Xi(J))$ has a trindex for
a.e.~$\lambda\in\Lambda$. Then the {\it generalized spectral shift
function} $\what\xi(\cdot,H_0,H)$
associated with the pair $(H_0,H)$ is defined by
\begin{equation}
\what\xi(\lambda,H_0,
H)=\trind(\what\Xi(\lambda,H_0,H), \Xi(J)) \text{ for a.e. }
\lambda\in \Lambda. \lb{5.73}
\end{equation}
\end{definition}

\begin{remark} \lb{r5.9a}
A close look at $\phi(z)$ and $\what\Xi$ in \eqref{5.71} and \eqref{5.72}
reveals that both objects depend on the self-adjoint operators $H_0$ and
$V$. Thus, a logical choice of notation for $\what\Xi$ would have
indicated its dependence on the pair $(H_0,V)$. We decided against that
choice since in practical applications, $\what\Xi$ in \eqref{5.71} is
associated with a pair of self-adjoint operators $(H_0,H)$, where $H$
results as an additive perturbation of $H_0$ by $V$ and hence resorted
to the more familiar notation $\what\Xi(\lambda,H_0,H)$. But this raises
the question of how to define such a self-adjoint operator $H$, given
$H_0$ and $V$. Perhaps the most natural solution of this problem in our
context goes back to Kato \cite{Ka66}  (see also \cite{KK66}) and proceeds
as follows. One defines the resolvent of the self-adjoint operator $H$ in
$\calH$ (and hence $H$ itself) by
\begin{align}
(H-z)^{-1}&=(H_0-z)^{-1} \lb{5.73c} \\
& \quad -(|V|^{1/2}(H_0-\ol z)^{-1})^*
\phi(z)^{-1}\ol {|V|^{1/2}(H_0-z)^{-1}}, \quad 
 z\in\bbC_+. \no
\end{align}
A detailed discussion of this point of view can be found in Yafaev's
monograph \cite[Sects.~1.9, 1.10]{Ya92}.
\end{remark}

\vspace*{3mm}

By Theorem~\ref{t5.8}, the generalized spectral shift function coincides
with Krein's spectral shift function in the case of trace class
perturbations, that is,
\begin{equation}
\what\xi(\lambda,H_0,H)=\xi(\lambda,H_0,H) \text{ for a.e. }
\lambda\in\bbR, \lb{5.73a}
\end{equation}
using the standard factorization of $(H-H_0)$ into 
$(H-H_0)=|V|^{1/2}\sgn(V)|V|^{1/2}\in\calB_1(\calH)$.

\begin{lemma}\lb{l5.10}
Let $S$ be a signature operator, $S=S^*=S^{-1}$,
 $A=A^*\in \calB_{\infty}(\calH)$, and  $\Lambda=\bbR\backslash
\{\spec(S+A)\cup \{-1,1\}\}$. Then the generalized spectral
shift function $\what\xi (\lambda, S+A,S)$ associated with 
the pair $(S+A,S)$ is well-defined for a.e.
$\lambda\in\Lambda$. Moreover, $\what\xi (\lambda,
S+A,S)$ has a continuous representative on $\Lambda$ \rm{(}still denoted
by
$\what\xi (\lambda,S+A,S)$\rm{)} such that
\begin{equation}
\what\xi (\lambda, S+A,S)=\ind(E_{S+A}((-\infty,\lambda )),
E_{S}((-\infty,\lambda)) ),\quad \lambda\in \Lambda.
\lb{5.74}
\end{equation}
In particular, taking $\lambda\uparrow 0$,
\begin{equation}
\what\xi (0_-, S+A,S)=\ind(E_{S+A}((-\infty,0 )), E_{S}((-\infty, 0)) ),
\quad \lambda\in \Lambda.
\lb{5.75}
\end{equation}
\end{lemma}
\begin{proof}
Since $\spec(S)\subseteq \{-1,1\}$ and $A\in\calB_\infty(\calH)$, the
spectrum  of $S+A$ is a discrete set with only possible accumulation
points at
$\pm 1$. Next, the normal boundary values
$|A|^{1/2}(S+A-\lambda+i0)^{-1}|A|^{1/2}=
|A|^{1/2}(S+A-\lambda)^{-1}|A|^{1/2}$ exist in norm for all
$\lambda\in\Lambda$. Moreover,
\begin{equation}
(\sgn(-A)+|A|^{1/2}(S+A-\lambda)^{-1}|A|^{1/2})^{-1}\in\calB(\calH),
\quad \lambda\in\Lambda, \lb{5.76}
\end{equation}
which can be seen as follows: suppose that \eqref{5.76} is false,
then by compactness of $A$ there exists an $f\in\calH$ such that
\begin{equation}
(\sgn(-A)+|A|^{1/2}(S+A-\lambda)^{-1}|A|^{1/2})f=0. \lb{5.76a}
\end{equation}
Multiplying \eqref{5.76a} by $\sgn(-A)$ one infers that
\begin{equation}
(I_\calH+\sgn(-A)|A|^{1/2}(S+A-\lambda)^{-1}|A|^{1/2})f=0.
\lb{5.76b}
\end{equation}
Since
$A\in\calB_\infty(\calH)$ and
$\spec(CD)\backslash\{0\}=\spec(DC)\backslash\{0\}$ for any
$C,D\in\calB(\calH)$, one concludes that there is a $g\in\calH$ such that
$(I_\calH -(S+A-\lambda)^{-1}A)g=0$. Thus, $(S-\lambda)g=0$ and hence
$\lambda\in\{-1,1\}$,  which
contradicts the fact that $\lambda\in\Lambda$. This proves \eqref{5.76}.

By Theorem~\ref{t5.8a},
\begin{equation}
\what\Xi(\lambda,S+A,S)=
\Xi(\sgn(-A)+|A|^{1/2}(S+A-\lambda)^{-1}|A|^{1/2})
\quad \text{ for a.e. } \lambda\in\Lambda. \lb{5.78}
\end{equation}
Moreover, the pair
\begin{equation}
(\sgn(-A)+|A|^{1/2}(S+A-\lambda)^{-1}|A|^{1/2}, \sgn(-A)) \lb{5.77}
\end{equation}
is a Fedholm pair. Hence the generalized spectral shift function is
well-defined and given by
\begin{align}
&\what\xi(\lambda,
S+A,S)=\trind(\Xi(\sgn(-A)+|A|^{1/2}(S+A-
\lambda)^{-1}|A|^{1/2}),\Xi(-A)) \no \\
&=\ind(E_{\sgn
(-A)+|A|^{1/2}(S+A-\lambda)^{-1}|A|^{1/2}}((-\infty, 0)),
E_{\sgn(-A)}((-\infty, 0))) \lb{5.79} \\
& \hspace*{8.5cm}\text{ for a.e. } \lambda\in\Lambda. \no
\end{align}
Since the right-hand side of \eqref{5.79} is continuous on $\Lambda$
by Theorem~\ref{ttr.7}, $\what\xi (\lambda,S+A,S)$ has a continuous
representative on $\Lambda$.

Next, assume  $A\in \calB_1(\calH)$. Then Krein's spectral shift function
$\xi(\lambda,S+A,S)$ associated with the pair $(S+A,S)$ coincides with
the right-hand side of \eqref{5.74} for a.e.~$\lambda\in\bbR$ (see, e.g.,
\cite{BP98}), proving \eqref{5.74} for $A\in \calB_1(\calH)$ applying
\eqref{5.73a}.

The general case of compact operators $A\in\calB_\infty(\calH)$ can be
handled using an appropriate approximation argument. Denoting by
$\{\lambda_n\}_{n\in \bbZ}$ the eigenvalues of $A$ and by
$\{P_n\}_{n\in \bbZ}$
 the corresponding spectral  projections associated with
$\lambda_n$, and introducing the family of the self-adjoint operators
\begin{equation}
A_\rho=\sum_{n\in \bbZ}\rho^{-|n|}\lambda_nP_n, \quad \rho\in (0,1),
\lb{5.80}
\end{equation}
one concludes that $A_\rho\in \calB_1(\calH)$, $\rho\in (0,1)$, and
\begin{equation}
\nlim_{\rho\uparrow 1} \|A_\rho-A\|=0. \lb{5.81}
\end{equation}
Given $\lambda\in \Lambda$, there exists a $\rho_0\in (0,1)$, such that
for all $\rho\in (\rho_0, 1)$ the point
$\lambda\in\Lambda_\rho$, $\Lambda_\rho=\bbR\backslash\{\spec
(S+A_\rho)\cup \{-1,1\}\}$, and therefore, by \eqref{5.74} (for
$A\in\calB_1(\calH)$), for such $\rho$ we have the representation
\begin{equation}
\xi (\lambda, S+A_\rho,S)
=\ind(E_{S+A_\rho}((-\infty,\lambda )), E_{S}((-\infty, \lambda)) ),
\quad \rho\in (\rho_0,1). \lb{5.82}
\end{equation}
Here $\xi (\lambda, S+A_\rho,S)$ denotes the continuous representative
of (the piecewise constant) Krein's spectral shift function on
$\Lambda_\rho$.

Applying Theorem~\ref{ttr.7} once again,
one can pass to the limit $\rho \uparrow 1$ to obtain
\begin{equation}
\lim_{\rho\uparrow 1}\xi (\lambda, S+A_\rho,S)=
\ind(E_{S+A}((-\infty,\lambda )), E_{S}((-\infty, \lambda)) ,
\quad \lambda\in \Lambda. \lb{5.83}
\end{equation}
By \eqref{5.79} we also have
\begin{align}
&\xi(\lambda, S+A_\rho,S)
\no \\
&=\ind(E_{\sgn (-A_\rho)+|A_\rho|^{1/2}(S+A_\rho-
\lambda)^{-1}|A_\rho|^{1/2}}((-\infty, 0)),
E_{\sgn(-A_\rho)}((-\infty, 0)) )\lb{5.84}.
\end{align}
Taking into account that
$ \sgn(-A_\rho)=\sgn(-A)$, $\rho\in (0, 1)$ \eqref{5.84}
implies
\begin{align}
&
\lim_{\rho\uparrow 1} \xi(\lambda, S+A_\rho,S)
\no \\
&=
\ind(E_{\sgn (-A)+|A|^{1/2}(S+A-
\lambda)^{-1}|A|^{1/2}}((-\infty, 0)),
E_{\sgn(-A)}((-\infty, 0)) ),\quad \lambda\in \Lambda \lb{5.85}
\end{align}
by Theorem~\ref{ttr.7}. The right-hand side of \eqref{5.85}
coincides with the continuous representative of the generalized spectral
shift function $\what\xi(\lambda, S+A,S)$, $\lambda\in\Lambda$, which
together with \eqref{5.83} proves \eqref{5.74}. Finally, \eqref{5.75}
is a consequence of \eqref{5.74} and the left continuity of
$\ind(E_{S+A}((-\infty,\lambda )), E_{S}((-\infty, \lambda)) )$ on
$\bbR\setminus \{-1,1\}$.
\end{proof}

Combining Theorem~\ref{main} and Lemma~\ref{l5.10}, one can
finally reformulate Theorem~\ref{main} as follows
using the concept of the generalized spectral shift function.

\begin{theorem}\lb{BSCHW}
Under the assumptions of Theorem~\ref{main},
the spectral shift function $\what\xi(\lambda, H_0,H)$ associated 
with the pair $(H_0,H)$ admits the representation
\begin{equation}\lb{bsc}
\what\xi(\lambda, H_0,H)= \frac{1}{\pi}\int_\bbR dt\,
\frac{\what\xi(0_-,J+A(\lambda)+tB(\lambda),J)}{1+t^2}, \text{ for
a.e. } \lambda \in \bbR,
\end{equation}
where $\what\xi(\cdot,J+A(\lambda)+tB(\lambda),J)$ is the continuous
representative of the generalized spectral shift function
associated with the pair
$ (J+A(\lambda)+tB(\lambda),J)$ for a.e. $\lambda\in\bbR$, $ t\in \bbR$.
\end{theorem}
\begin{proof}
The assertion is a direct consequence of Theorem~\ref{main} and
Lemma~\ref{l5.10}.
\end{proof}

\begin{remark}\lb{principle}
Suppose $H, H_0, V$ are self-adjoint in $\calH$, with
$V\in\calB_1(\calH)$ and $H=H_0+V$. Suppose $H_0$ has a spectral
gap and $\lambda\in \Lambda$, with $\Lambda$ a joint spectral gap of
$H_0$ and $H$. Then
\eqref{bsc} turns into
\begin{equation}\lb{bsn}
\xi(\lambda, H_0,H)=
\xi(0,J+|V|^{1/2}(H_0-\lambda)^{-1}|V|^{1/2},J), \quad
\lambda\in\Lambda,
\end{equation}
where $\xi(\lambda,H_0,H)$
($\xi(0,J+|V|^{1/2}(H_0-\lambda)^{-1}|V|^{1/2},J)$) denotes the
continuous representative of  Krein's spectral shift function associated
with the pair $(H_0,H)$ ($(J+|V|^{1/2}(H_0-\lambda)^{-1}|V|^{1/2},J)$)
on $\Lambda$. In particular, if $H_0$ is bounded from below and
$\lambda<\inf(H_0)$, and the perturbation $V$ is non-positive (i.e.,
$V\le0)$, the equality \eqref{bsn} has the following meaning: the number
of  eigenvalues of the operator $H=H_0+V$, located to the left of the
point $\lambda $, $\lambda<\inf(H_0)$, coincides with the number
of eigenvalues of $|V|^{1/2}(H_0-\lambda)^{-1}|V|^{1/2}$ which are greater
than $1$. Therefore, in this special case where
$\lambda<\inf\spec(H_0)$, \eqref{bsn} represents the classical
 {\it Birman--Schwinger principle} (a term coined by Simon, see, e.g.,
\cite{Si79})  as originally introduced by  Birman
\cite{Bi66} (see also \cite[Ch.~7]{Sc61}). Thus, \eqref{bsn} should be
interpreted as the {\it Birman-Schwinger principle in a gap} and hence
\eqref{bsc} as the {\it  generalized Birman-Schwinger principle.} We
emphasize that Theorem~\ref{BSCHW} introduces a new twist in connection
with   the  (generalized) Birman-Schwinger principle: the role of
eigenvalue counting
functions in the traditional formulation of the
Birman-Schwinger principle (see \cite{BS91}
for a modern formulation of the principle) is now replaced by the more
general concept of  the (generalized)  spectral shift function
$\what\xi(\lambda, H_0,H)$ and an appropriate 
average over $\what\xi(0_-,J+A(\lambda)+tB(\lambda),J)$.
\end{remark}


\end{document}